\documentclass[a4paper,11pt]{amsart}
\usepackage[utf8]{inputenc}

\usepackage{lipsum}

\usepackage{enumerate, mathrsfs}
\usepackage[dvipsnames,usenames]{xcolor}

\usepackage{enumitem}

\usepackage{amsmath,amsfonts,amsthm,amssymb}
\usepackage[all]{xy}
\usepackage[T1]{fontenc}

\makeatletter
\newcommand{\addresseshere}{%
  \enddoc@text\let\enddoc@text\relax
}
\makeatother

\usepackage{hyperref}
\hypersetup{colorlinks=true, urlcolor=Blue, citecolor=Green, linkcolor=BrickRed, breaklinks, unicode}

\usepackage{cleveref}
\usepackage{tikz}

\usepackage{verbatim}

\usepackage{booktabs}

\usepackage[font=small,labelfont=bf]{caption}

\theoremstyle{plain}

\theoremstyle{definition}

\DeclareMathOperator{\crs}{cr}
\DeclareMathOperator{\maxcr}{Top}
\DeclareMathOperator{\extra}{m}
\DeclareMathOperator{\recal}{R}

\title{Hard Diagrams of the Unknot}

\date{}

\author[B]{Benjamin A. Burton}
\address{Benjamin A. Burton: School of Mathematics and Physics, The University of Queensland, St Lucia QLD 4072, Australia}\email{bab@maths.uq.edu.au}

\author[C]{Hsien-Chih Chang}
\address{Hsien-Chih Chang: Department of Computer Science, Dartmouth College, Hanover, NH 03755, United States}
\email{hsien-chih.chang@dartmouth.edu}

\author[L]{Maarten L\"offler}
\address{Maarten L\"offler: Department of Information and Computing Sciences, Utrecht University, Buys Ballot Gebouw, 
Princetonplein, 53584 CC Utrecht, The Netherlands}
\email{m.loffler@uu.nl}

\author[deM]{Arnaud de Mesmay}
\address{Arnaud de Mesmay: Laboratoire d'Informatique Gaspard Monge, 5 Boulevard Descartes, Champs-sur-Marne, 77454 Marne-la-Vallée Cedex 2, France}
\email{arnaud.de-mesmay@univ-eiffel.fr}

\author[M]{Cl\'ement Maria}
\address{Cl\'ement Maria: DataShape Research Group, INRIA Sophia Antipolis-M\'editerran\'ee, 2004 route des Lucioles, 06902 Sophia Antipolis, France}
\email{clement.maria@inria.fr}

\author[S]{Saul Schleimer}
\address{Saul Schleimer: Mathematics Institute, University of Warwick, Coventry CV4 7AL, United Kingdom}
\email{s.schleimer@warwick.ac.uk}

\author[S]{Eric Sedgwick}
\address{Eric Sedgwick: School of Computing, DePaul University, 2400 N Sheffield Ave, Chicago, IL 60614, United States}
\email{ESedgwick@cdm.depaul.edu}

\author[S]{Jonathan Spreer}
\address{Jonathan Spreer: School of Mathematics and Statistics F07, The University of Sydney, NSW 2006, Australia} \email{jonathan.spreer@sydney.edu.au}

\begin{document}

\maketitle

\begin{abstract}
  We present three ``hard'' diagrams of the unknot. They require (at least) three extra crossings before they can be simplified to the trivial unknot diagram via Reidemeister moves in~$\mathbb{S}^2$. Both examples are constructed by applying previously proposed methods. The proof of their hardness uses significant computational resources.
  We also determine that no small ``standard'' example of a hard unknot diagram requires more than one extra crossing for Reidemeister moves in~$\mathbb{S}^2$.
\end{abstract}

\section{Introduction}

Hard diagrams of the unknot are closely connected to the complexity of the unknot recognition problem. Indeed a natural approach to solve instances of the unknot recognition problem is to untangle a given unknot diagram by Reidemeister moves. Complicated diagrams may require an exhaustive search in the Reidemeister graph that very quickly becomes infeasible.

As a result, such diagrams are the topic of numerous publications~\cite{Freedman1994energy,Goeritz34,Henrich14,kauffman2012hard,Ochiai90}, as well as discussions by and resources provided by leading researchers in low-dimensional topology~\cite{AgolCollection,CurtisCollection,MathOverflow_1,MathOverflow_2}.
However, most classical examples of hard diagrams of the unknot are, in fact, easy to handle in practice. Moreover, constructions of more difficult examples often do not come with a rigorous proof of their hardness~\cite{Freedman1994energy,kauffman2012hard}.

The purpose of this note is hence two-fold. Firstly, we provide three unknot diagrams together with a proof that they are significantly more difficult to untangle than classical examples. Secondly, we survey existing literature on hard unknots and compare the hardness of well-known examples, thereby making claims made in, for instance, \cite[Pages 41--42]{Freedman1994energy} and \cite[Page 4--5 and Section 9]{kauffman2012hard} more comparable.

We want this work to be practically usable. For this reason, we provide Gauss codes for each of the unknot diagrams discussed in this note in~\Cref{appendix}. These codes should be suitable as input for most software on knots with at most minimal adjustments. Their current format is compatible with \emph{Regina}~\cite{regina}.

\medskip

If you are familiar with the topic of finding difficult diagrams of the unknot, and you are mainly interested in our examples $D_{28}$, $D_{43}$, and $PZ_{78}$, please skip ahead to \Cref{sec:hard} or \Cref{appendix}.

\subsection{Basics about knots}
\label{sec:basics}

\begin{figure}
  \centerline{\includegraphics[width=0.8\textwidth]{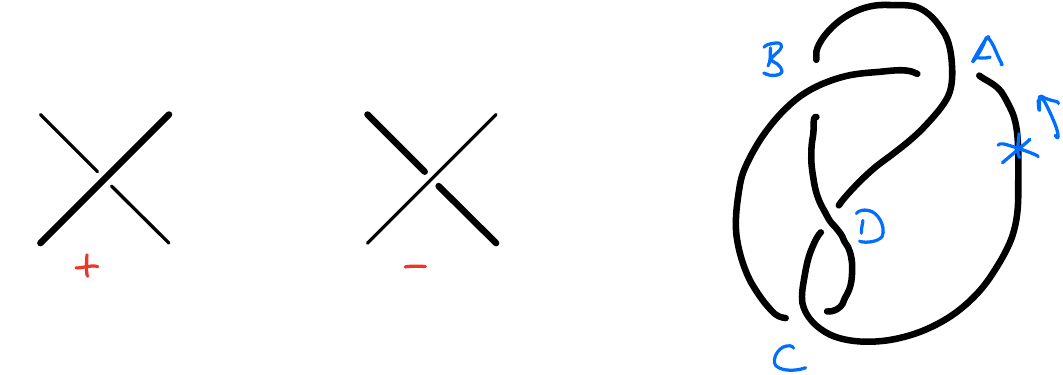}}
  \caption{Left: signs of a crossing. Right: the $4$-crossing figure-eight knot diagram. Its Gauss code with respect to $\star$ is
  $(-A)(+B)(-C)(+D)(-B)(+A)(-D)(+C)$.}
\label{fig:Gauss}
\end{figure}

A \emph{knot} is a piecewise-linear embedding of $S^1$ into the three-dimensional sphere $\mathbb{S}^3$. A given knot $K \subset \mathbb{S}^3$ is often represented by a {\em diagram} $D_K$, that is, a projection of $K$ into the $2$-sphere $\mathbb{S}^2$.%
\footnote{Note that this definition of a knot diagram may not be considered as standard. In many settings the target of the projection is the plane $\mathbb{R}^2$ rather than its one-point compactification $\mathbb{R}^2 \cup \{\infty\} = \mathbb{S}^2$; this is more intuitive when drawing a diagram by hand. We focus on knot diagrams in the $2$-sphere because it is more natural from an algorithmic point of view.}
In such a knot diagram, the projection must be in general position admitting at most a finite number of double points, called {\em crossings}. Their number is denoted by $\crs(D_k)$. At every crossing the diagram provides additional information on which local segment of the knot goes ``over'' the other.
In addition, a knot diagram is usually endowed with an {\em orientation}, that is, a choice of direction in which we run along the knot, indicated by an arrow. However, since a choice of orientation has no effect on the hardness of a knot diagram, diagrams in this note are non-oriented.
See \Cref{fig:Gauss,fig:hard_example,fig:harder_example} for examples of knot diagrams.

Knot diagrams can be represented purely combinatorially in various ways. One of them is what is called the {\em Gauss code}: Giving a knot diagram with labels at the crossings we run along the knot writing down the labels of the crossings as we encounter them together with a sign encoding whether we go over ($+$) or under ($-$).
By construction, we see every crossing twice with opposite signs. The Gauss code determines a knot type up to orientation. This is sufficient for our purposes. See \Cref{fig:Gauss} for an example.

\begin{figure}[h]
\centerline{\includegraphics[width=\textwidth]{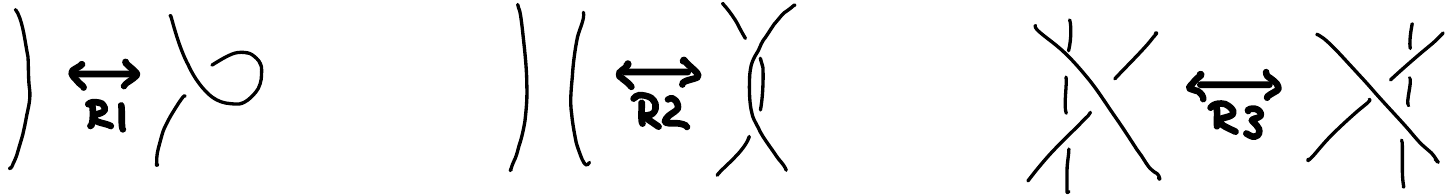}}
\caption{The three types of Reidemeister moves R1, R2 and R3 together with their inverses. \label{fig:R1R2R3}}
\end{figure}

Two knots $K_1$ and $K_2$ are considered \emph{equivalent} if there is a continuous deformation (an \emph{isotopy}) between them respecting the embeddings of $K_1$ and $K_2$. It follows from a theorem by Reidemeister~\cite{Reidemeister27} that two knot diagrams represent equivalent knots if and only if one can be transformed into the other by a sequence of local modifications of the diagrams called {\em Reidemeister moves}.
See \Cref{fig:R1R2R3} for an illustration of all three Reidemeister moves and their inverses, and see the first two steps in \Cref{fig:culprit_undone} for a version of R1 in~$\mathbb{S}^2$ (this can be seen as taking an outermost arc of $D$, and dragging it over the entire diagram by passing through $\infty$).

Due to Reidemeister's theorem we can also define the {\em Reidemeister graph} of a knot, where nodes are $\mathbb{S}^2$-isotopy classes of diagrams of a knot and two nodes are connected by an arc if and only if their corresponding diagrams are transformed into each other by a single Reidemeister move.

The unique knot admitting a diagram with no crossings is called the {\em unknot}. It's trivial $0$-crossing diagram is denoted by $D_{0}$.
Let $D$ be a diagram of the unknot with $\crs(D)$ crossings. By Reidemeister's theorem, $D$ admits a sequence of Reidemeister moves connecting it to the trivial diagram $D_{0}$. Naturally, in every such sequence there is an intermediate diagram with the largest number of crossings. For typical input with few crossings this will often be the initial diagram itself.

Fixing $D$, we are interested in the sequence of Reidemeister moves minimising this maximal crossing number. In accordance with Kauffman and Lambropoulou~\cite[Page 5]{kauffman2012hard} we denote this maximal crossing number by $\maxcr(D)$. Naturally one has $\maxcr(D) \ge \crs(D)$. We then say that $D$ requires $\extra(D) = \maxcr(D) - \crs(D)$ {\em extra crossings} for it to be simplified to $D_{0}$. In this note we use $\extra (D)$ as a measure of hardness for unknot diagrams.
(Note that Kauffman and Lambropoulou~\cite{kauffman2012hard} define a different hardness measure for unknot diagrams by $\recal(D) = \maxcr(D) / \crs(D)$, which they call {\em recalcitrance}.)

Naturally, the notion of extra crossings can be extended to diagrams of arbitrary knots where the task is then to simplify a given diagram to one with the minimum number of crossings.

\subsection{Background on difficult unknot diagrams}
\label{sec:background}

This note is about diagrams of the unknot. Specifically, we are interested in the following question:

\medskip

\begin{center}
\textbf{Question 1}: {\em Are there diagrams of the unknot $D$ such that $\extra (D) \gg 0$?}
\end{center}

\medskip

This question is fundamental in algorithmic knot theory and has attracted lots of interest: simplifying an input knot diagram comes as a very first step in virtually every computation on knots. The increase in the size of intermediate diagrams, i.e., the number of extra crossings, measures the difficulty of the search for a diagram with fewer crossings.
The minimal number of Reidemeister moves to untangle an unknot diagram, and the maximal crossing number of any intermediate diagram met along the way, are related.

Early exponential upper bounds on the number of Reidemeister moves to untangle an unknot diagram~\cite{Hass01Bound} proved that the number of extra crossings is bounded by $\extra(D) \leq 2^{c \cdot \crs(D)}$, for a large constant $c$.
Subsequent work by Dynnikov~\cite{Dynnikov06ArcPresentation} on arc presentations implied the existence of a super-polynomially long sequence of Reidemeister moves on an unknot diagram $D$, leading to the trivial diagram, such that the maximal number of crossings of intermediate diagrams is bounded by a quadratic function:
\[
  \maxcr(D) = \crs(D) + \extra(D) \leq (\crs(D) - 1)^2 /2.
\]
More recently, Lackenby proved in~\cite{Lackenby13PolyBoundReidemeister} that unknot diagrams can be simplified with only polynomially many Reidemeister moves (more precisely, $O(\crs(D)^{11})$) without ever exceeding a quadratic number of crossings; in other words, $\maxcr(D) = O(\crs(D)^2)$.

Regarding lower bounds, Hass and Nowik exhibit an infinite family of unknot diagrams that require at least a quadratic number of Reidemeister moves to be untangled~\cite{Hass08BoundBetweenDiagrams,Hass10QuadraticBound}. However, these diagrams can be untangled with a monotonically non-increasing number of crossings, hence $\extra(D) = 0$.
In fact, the ``hardest'' unknot diagrams discussed in the literature admit straightforward simplifying sequences of Reidemeister moves that create only a few extra crossings --- if any at all. Moreover, given such a ``hard'' unknot diagram, the exact number of Reidemeister moves required is only provided in a small number of cases:
For Reidemeister moves in the plane (and not in $\mathbb{S}^2$, the setup of this article), there are claims of a $10$-crossing unknot diagram requiring two extra crossings~\cite{kauffman2012hard} (known as the {\em Culprit}), and a $32$-crossing unknot diagram believed to require four extra crossings~\cite[Pages 41--42]{Freedman1994energy} (which we call the \textit{Freedman-He-Wang} unknot).
Also, an infinite family of unknot diagrams is claimed to have unbounded recalcitrance (and thus super-linearly many extra crossings)~\cite{kauffman2012hard}.

\subsection{Our Contributions}
\label{sec:contributions}

Our contributions to Question 1 are the following:

\smallskip

In \Cref{sec:errors} we show that several hard diagrams of the unknot are, in fact, not very hard, when Reidemeister moves in $\mathbb{S}^2$ are considered. More precisely, for most hard unknot diagrams in the literature, at most one extra crossing is needed for them to be untangled (that is, we have $\extra (D) \leq 1$ in these cases). Specifically, both the Culprit from Kauffman and Lambropoulou~\cite{kauffman2012hard}
and the Freedman-He-Wang unknot~\cite{Freedman1994energy} only require one extra crossing.
Moreover, we raise concerns about the completeness of the proof and the correctness of the statement about the family of unknot diagrams in~\cite{kauffman2012hard} requiring a super-linear number of extra crossings. See \Cref{ssec:inaccuracies} for details.

\smallskip

In \Cref{sec:hard} we present three diagrams $D_{28}$, $D_{43}$, and $PZ_{78}$ of the unknot, and give a computer proof that $\extra(D_{28}) = 3$, $\extra (D_{43}) \geq 3$, and $\extra (PZ_{78}) \geq 3$ even for Reidemeister moves in $\mathbb{S}^2$.

\smallskip

Finally, in \Cref{sec:harder} we present the approaches that led to the discovery of two of our hard unknot diagrams, $D_{28}$ and $D_{43}$\footnote{Our third hard unknot diagram $PZ_{78}$ was discovered by simplifying the unknot diagram from~\cite[Figure~27]{petronio2016algorithmic} and is not discussed in further detail in this article.}. We believe that both of these approaches are suitable to construct infinite families of unknot diagrams requiring an unbounded number of extra crossings.

\subsection*{Acknowledgements}

The authors would like to thank Dagstuhl and the organizers of the Dagstuhl seminar 19352: ``Computation in Low-Dimensional Geometry and {Topology}'' where this work has been initiated.

\section{Current literature}
\label{sec:errors}

\subsection{Difficult unknots in the literature}
\label{ssec:literature}

Many examples of diagrams of the unknot that are not straightforward to untangle can be found in the literature~\cite{Freedman1994energy,Goeritz34,Henrich14,kauffman2012hard,Ochiai90,petronio2016algorithmic}.
In this note we focus on the number of extra crossings $\extra (D)$ necessary to untangle an unknot diagram $D$ (using Reidemeister moves in $\mathbb{S}^2$, see \Cref{sec:basics} for a precise definition) as a measure of hardness, and assess well-known examples against this measure.
Note that, while $\extra(D) > 0$ is relatively easy to establish, $\extra (D) \geq 2$ is much harder to verify, as the portion of the Reidemeister graph required to exhaustively search through expands rapidly in size.

Find below a list of hard unknot diagrams $D$, their number of crossings, their number of extra crossings $\extra(D)$, and where they can be found in the literature.  Indefinite answers (e.g., $\extra(D) \geq 2$) mean that the portion of the Reidemeister graph necessary to exhaustively search through does not fit into $\sim$8GB of memory.
(See \Cref{appendix} for Gauss codes of each of these knots.)

\medskip

\begin{center}
\begin{tabular}{llll}
  \toprule
  Name & $\crs$ & $\extra(D)$ & References \\
  \midrule
  H	 		&$9$   &$1$ 	&	\cite[Figure~4]{Henrich14} \\
  J	 		&$9$   &$1$ 	&	\cite[Figure~4]{Henrich14} \\
  Culprit 		&$10$  &$1$	&	\cite{kauffman2012hard} \\
  Monster 		&$10$  &$0$	&	 \cite{petronio2016algorithmic, CurtisCollection}	\\
  Goeritz 		&$11$  &$1$	&	\cite{Goeritz34}, cf.\ \S\ref{ssec:goeritz} \\
  Thistlethwaite 	&$15$  &$0$	&	\cite[Figure~9]{petronio2016algorithmic}, \cite{CurtisCollection}\\
  Ochiai I 		&$16$  &$0$	&	\cite[Figure 1]{Ochiai90} \\
  Tuzun-Sikora		&$21$  &$0$	&	\cite[Figure 8]{Tuzun16} \\
  Freedman-He-Wang 		&$32$  &$0$	&	\cite[Figure 6.1]{Freedman1994energy}, \S\ref{ssec:freedman} \\
  ``Fake'' Freedman-He-Wang 		&$32$  &$0$	&	\S\ref{ssec:freedman} \\
  Ochiai II 		&$45$  &$0$	&	\cite[Figure~2]{Ochiai90} \\
  Ochiai II (reduced) 	&$35$  &$\geq 2$&	\S\ref{sec:errors} \\
  Ochiai III 	&$67 $  &$ 0$&		\cite[Figure~3]{Ochiai90} \\
  Ochiai IV (Suzuki) 		&$55$  &$\geq 2$&	\cite[Figure~4]{Ochiai90} \\
  Haken 		&$141$ &$0$	&	\cite{petronio2016algorithmic} \\
  $PZ_{31}$ 	&$31$ &$0$	&	\cite[Figure~12]{petronio2016algorithmic} \\
  $PZ_{120}$	&$120$ &$\geq 2$	&	\cite[Figure~27]{petronio2016algorithmic} \\
  $PZ_{138}$	&$138$ &$\geq 2$	&	\cite[Figure~14]{petronio2016algorithmic} \\
  \bottomrule
\end{tabular}
\end{center}

By Ochiai I, II, III, and IV we mean the unknot diagrams given in Figures 1, 2, 3, and 4 in~\cite{Ochiai90} respectively. According to the author, the diagram given in Figure~4 of~\cite{Ochiai90} (Ochiai IV) is a slightly modified version of a diagram originally due to Suzuki.

Note that the Monster, Thistlethwaite, Ochiai I, Tuzun-Sikora, Ochiai II, $PZ_{31}$, and Haken do not require extra crossings in order to be simplified. Some of them often reduce to smaller diagrams that then, in turn, require extra crossings to be simplified further:

Ochiai II reduces to a $35$-crossing diagram of the unknot for which at least two extra crossings are necessary for it to untangle. A slight modification of Freedman-He-Wang (referred to as the \textit{``Fake'' Freedman-He-Wang} unknot diagram, see \Cref{fig:fakeFHWD})  reduces to our $28$-crossing hard unknot diagram, denoted by $D_{28}$, see \Cref{sec:hard}, with $\extra(D_{28}) = 3$.
Moreover, Goeritz reduces to H (and, according to~\cite{kauffman2012hard}, also to the Culprit).

Example $PZ_{120}$ is too large for exhaustive enumeration. Example $PZ_{138}$ reduces to a $78$-crossing diagram $PZ_{78}$ using only two extra crossings, see \Cref{sec:hard}. The number $\extra (PZ_{78})$ is still unknown, and hence $\extra (PZ_{138}) \geq 2$. Exploring the Reidemeister graph around $PZ_{78}$ reveals more unknot diagrams requiring at least three extra crossings to be simplified. Is this a property of $PZ_{78}$ and its position in the Reidemeister graph, or is this common for unknot diagrams of around $78$ crossings? Answering this question in a meaningful may give new insight on the problem of recognising the unknot.

\subsection{Specific claims from the literature}
\label{ssec:inaccuracies}

In this section we comment on a number of claims found in the literature with respect to the hardness of some well-known unknot diagrams.

\paragraph{The Culprit:} In~\cite{kauffman2012hard} the authors make the following claim about simplifying sequences of Reidemeister moves (in $\mathbb{R}^2$) starting with the Culprit: ``\emph{...\ the largest increase is to a diagram of $12$ crossings. This is the best possible result for this diagram}''. They then go on to provide such a sequence.

\begin{figure}[h]
\centerline{\includegraphics[width=0.7\textwidth]{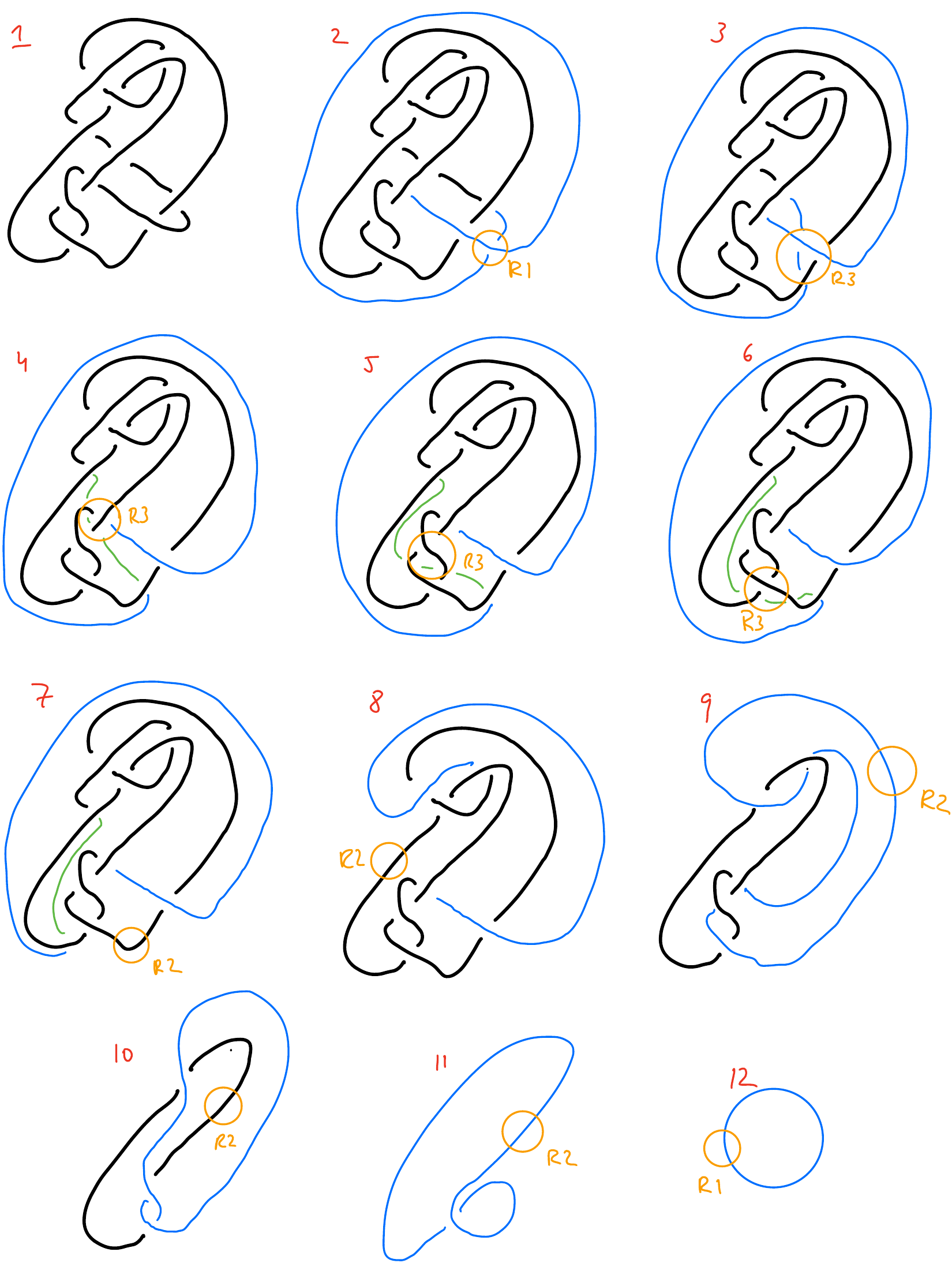}}
\caption{Untangling the Culprit with one extra crossing. \label{fig:culprit_undone}}
\end{figure}

When considering Reidemeister moves in the $2$-sphere, the Culprit can be transformed into the trivial unknot diagram with only one extra crossing; see \Cref{fig:culprit_undone}.

\paragraph{The Freedman-He-Wang unknot:} In~\cite{Freedman1994energy} the authors state that the Freedman-He-Wang unknot ``\emph{...\ cannot be connected to a round diagram by a family of knot diagrams without the number of crossings increasing to at least $33$}''. Even more, they conjecture that ``\emph{...\ from the picture it looks like a maximum of $36$ crossings must occur}'' in any unknotting sequence.

When considering Reidemeister moves in the $2$-sphere, there exists a sequence of Reidemeister moves turning Freedman's $32$-crossing unknot diagram into the trivial diagram by just using one extra crossing. Moreover, all intermediate diagrams require at most two extra crossings to be simplified further.

However, when considering a slight modification of the Freedman-He-Wang unknot (one that initially allows more R3 moves, which we refer to as the ``Fake'' Freedman-He-Wang unknot; see \Cref{fig:fakeFHWD}), we find a non-increasing sequence of Reidemeister moves down to our $28$-crossing unknot diagram $D_{28}$ (see \Cref{sec:hard}), which then, in turn, needs three extra crossings in order to be simplified to the trivial diagram.
We attribute the existence of this harder intermediate diagram to the fact that the ``Fake'' Freedman-He-Wang never expanded to $33$ crossings, while the original Freedman-He-Wang diagram did.

\paragraph{Infinite families of unbounded recalcitrance:} In~\cite{kauffman2012hard} the authors claim that there exist infinite families of type ``generalizations of Goeritz' unknot diagram'' with unbounded recalcitrance (or, equivalently, a number of extra crossings that is super-linear in the number of crossings of the initial diagram). This claim is never proven.
Instead, the authors argue that there exists an infinite family of generalizations of the Goeritz unknot such that the number of Reidemeister moves necessary to simplify these knots to the trivial diagram is quadratic in the initial number of crossings. In particular, it seems like the term recalcitrance is silently redefined to mean ``minimal number of Reidemeister moves necessary to simplify to $D_0$'' divided by ``number of crossings of diagram''.

However, note that examples of unknot diagrams requiring a quadratic number of Reidemeister moves in order to be simplified have already been presented in~\cite{Hass10QuadraticBound}.

\section{Three examples of ``hard'' unknot diagrams}
\label{sec:hard}

In this section we present three examples of hard unknot diagrams. More specifically, we present
\begin{itemize}
  \item a $28$-crossing diagram of the unknot, $D_{28}$, such that every sequence of Reidemeister moves untangling $D_{28}$ to the trivial unknot diagram $D_0$ requires three extra crossings (that is, we have $\extra (D_{28}) = 3$); and
  \item a $43$-crossing diagram of the unknot, $D_{43}$, requiring {\em at least} three extra crossings ($\extra (D_{43}) \geq 3$).
  \item a $78$-crossing diagram of the unknot, $PZ_{78}$, requiring {\em at least} three extra crossings ($\extra (PZ_{78}) \geq 3$).
\end{itemize}
The examples $D_{28}$ and $D_{43}$ were obtained by two distinct strategies. Both are briefly described in \Cref{sec:harder}. The example $PZ_{78}$ was obtained by simplifying a $138$-crossing unknot diagram due to Petronio and Zanelatti~\cite[Figure~27]{petronio2016algorithmic} by going through a $140$-crossing diagram.

The proof of the claimed properties of $D_{28}$, $D_{43}$, and $PZ_{78}$ goes by exhaustive enumeration in the Reidemeister graph. This is carried out using the knot interface of low-dimensional topology software {\em Regina}~\cite{regina}. The algorithm tests for all possible sequences of Reidemeister moves in the $2$-sphere.
Reproducing the results using {\em Regina} is straightforward. In the case of $D_{28}$ the required computing power is relatively moderate. This is different for $D_{43}$ and $PZ_{78}$.

The diagram $D_{28}$ is shown in \Cref{fig:hard_example}, the diagram $D_{43}$ is shown in \Cref{fig:harder_example}, and the diagram of $PZ_{78}$ is shown in \Cref{fig:PZ_example}, all in \Cref{appendix}. Their Gauss codes are given in the captions of their figures.

\section{Two methods to produce even ``harder'' unknot diagrams}
\label{sec:harder}

\subsection{Generalizing the Freedman-He-Wang unknot}
\label{ssec:freedman}

\begin{figure}[h]
      \centering
      \includegraphics[width=11cm]{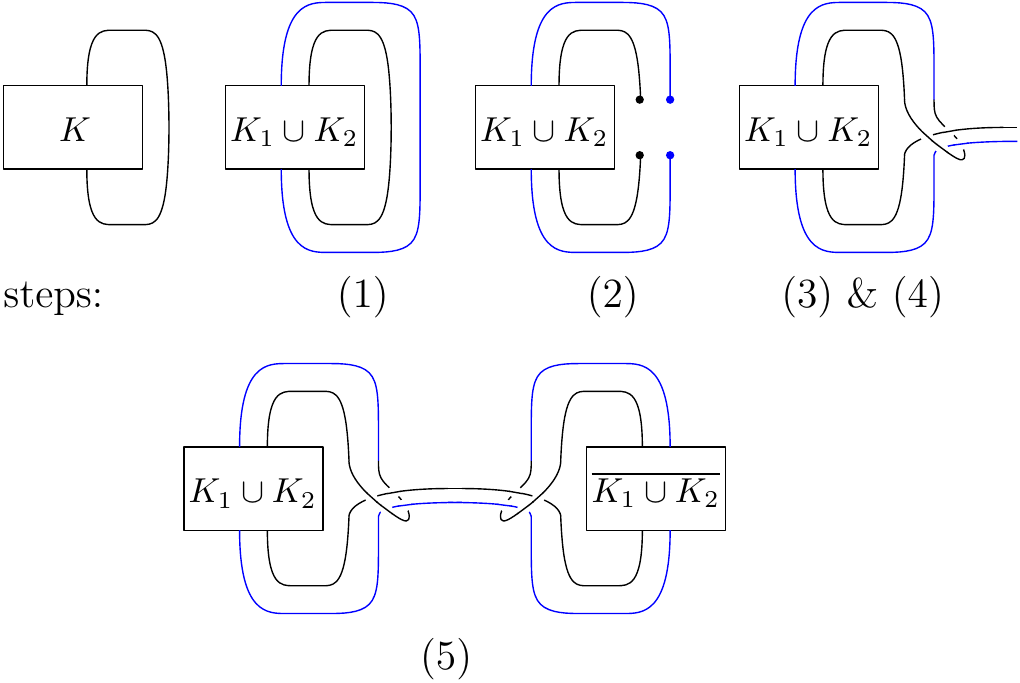}
      \caption{Construction generalizing the Freedman-He-Wang unknot diagram from~\cite{Freedman1994energy}.}
      \label{F:gen_Freedman}
 \end{figure}

\Cref{F:hard_unknot} pictures the classical hard unknot diagram due to Freedman, He, and Wang, which has been studied in the context of the energy minimization approach to the unknotting problem~\cite{Freedman1994energy}. It can be constructed and generalized with the following procedure. For a knot $K$,
\begin{enumerate}
\item    double the knot $K$ (along the blackboard framing) into two parallel copies $K_1$ and $K_2$,
\item    cut $K_1$ and $K_2$ open. We now have four ends coming in two pairs.
\item    Take two mirror symmetric pairs of ends and stretch them out as two parallel strands,
\item    take the remaining pair of ends and wrap them around those two strands before connecting them up,
\item    mirror a second copy of this tangle and build the connect-sum in the obvious way.
\end{enumerate}
This process is illustrated in \Cref{F:gen_Freedman}, and
\Cref{F:hard_unknot} shows the result of the construction for the trefoil knot $K$. This setup essentially led to the discovery of $D_{28}$, see \Cref{sec:hard}.

Proving lower bounds for generalizations of this example using knots with increasing bridge numbers is the object of ongoing work.

\subsection{Generalizing the Goeritz unknot}
\label{ssec:goeritz}

\begin{figure}[h]
      \centering
      \def\svgwidth{\textwidth}
      \begingroup%
  \makeatletter%
  \providecommand\color[2][]{%
    \errmessage{(Inkscape) Color is used for the text in Inkscape, but the package 'color.sty' is not loaded}%
    \renewcommand\color[2][]{}%
  }%
  \providecommand\transparent[1]{%
    \errmessage{(Inkscape) Transparency is used (non-zero) for the text in Inkscape, but the package 'transparent.sty' is not loaded}%
    \renewcommand\transparent[1]{}%
  }%
  \providecommand\rotatebox[2]{#2}%
  \newcommand*\fsize{\dimexpr\f@size pt\relax}%
  \newcommand*\lineheight[1]{\fontsize{\fsize}{#1\fsize}\selectfont}%
  \ifx\svgwidth\undefined%
    \setlength{\unitlength}{2375.32752496bp}%
    \ifx\svgscale\undefined%
      \relax%
    \else%
      \setlength{\unitlength}{\unitlength * \real{\svgscale}}%
    \fi%
  \else%
    \setlength{\unitlength}{\svgwidth}%
  \fi%
  \global\let\svgwidth\undefined%
  \global\let\svgscale\undefined%
  \makeatother%
  \begin{picture}(1,0.23243351)%
    \lineheight{1}%
    \setlength\tabcolsep{0pt}%
    \put(0,0){\includegraphics[width=\unitlength,page=1]{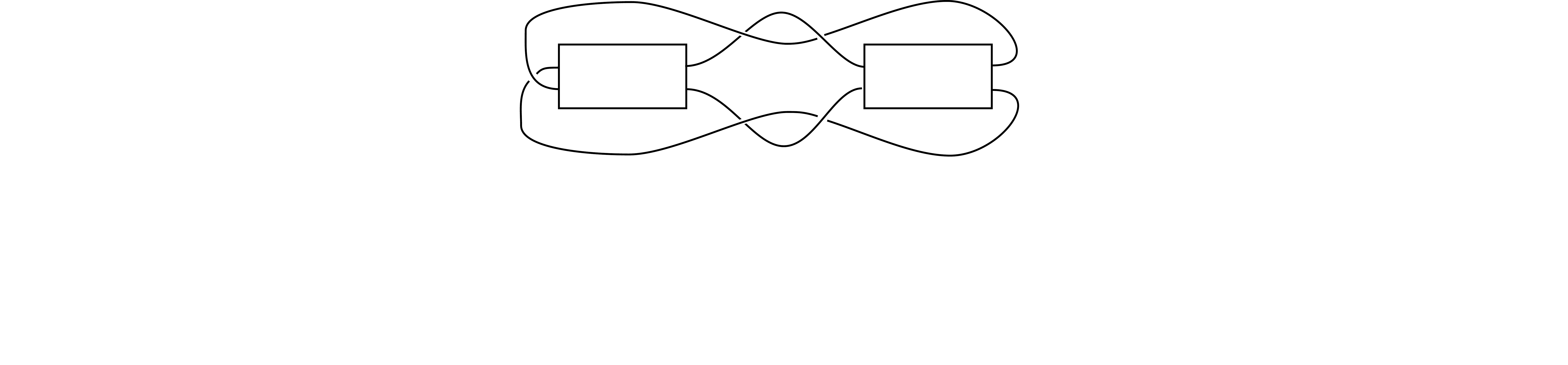}}%
    \put(0.39526885,0.17319216){\color[rgb]{0,0,0}\makebox(0,0)[t]{\lineheight{1.25}\smash{\begin{tabular}[t]{c}$B$\end{tabular}}}}%
    \put(0.59170792,0.17329392){\color[rgb]{0,0,0}\makebox(0,0)[t]{\lineheight{1.25}\smash{\begin{tabular}[t]{c}$B^{-1}$\end{tabular}}}}%
    \put(0,0){\includegraphics[width=\unitlength,page=2]{goeritz.pdf}}%
    \put(0.74529579,0.17319216){\color[rgb]{0,0,0}\makebox(0,0)[t]{\lineheight{1.25}\smash{\begin{tabular}[t]{c}\reflectbox{$B$}\end{tabular}}}}%
    \put(0.94173484,0.17329392){\color[rgb]{0,0,0}\makebox(0,0)[t]{\lineheight{1.25}\smash{\begin{tabular}[t]{c}$B^{-1}$\end{tabular}}}}%
    \put(0,0){\includegraphics[width=\unitlength,page=3]{goeritz.pdf}}%
    \put(0.4227188,0.04292727){\color[rgb]{0,0,0}\makebox(0,0)[t]{\lineheight{1.25}\smash{\begin{tabular}[t]{c}$B^{-1}$\end{tabular}}}}%
    \put(0.22447552,0.04282551){\color[rgb]{0,0,0}\makebox(0,0)[t]{\lineheight{1.25}\smash{\begin{tabular}[t]{c}$B$\end{tabular}}}}%
    \put(0,0){\includegraphics[width=\unitlength,page=4]{goeritz.pdf}}%
  \end{picture}%
\endgroup%
      \caption{The Goeritz unknot}
      \label{F:goeritz}
 \end{figure}

A classical example of a hard unknot due to Goeritz is pictured in \Cref{F:goeritz} on the top left. It can be thought of as the concatenation of two inverse braids with two flypes inserted in between on both of its strands (see the rest of \Cref{F:goeritz}). Undoing these flypes requires turning the braid which can be shown experimentally to require at least one additional crossing. 

\begin{figure}[h]
  \centering
  \includegraphics[width=5cm]{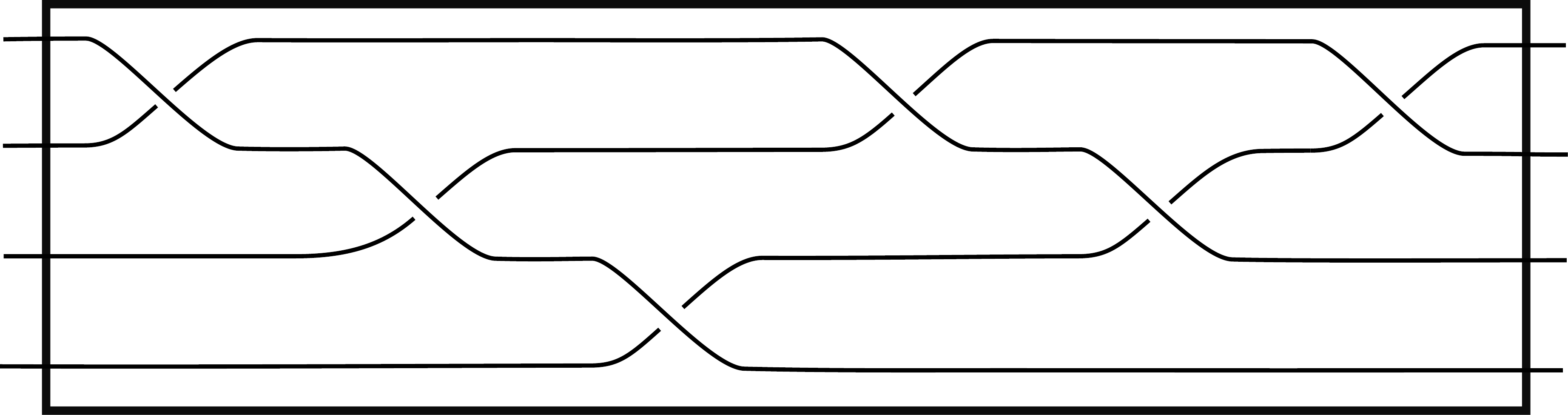}
  \caption{A four-strand generalization of a flype}
  \label{F:flype}
  \end{figure}

\begin{figure}[h]
      \centering
      \def\svgwidth{6cm}
\begingroup%
  \makeatletter%
  \providecommand\color[2][]{%
    \errmessage{(Inkscape) Color is used for the text in Inkscape, but the package 'color.sty' is not loaded}%
    \renewcommand\color[2][]{}%
  }%
  \providecommand\transparent[1]{%
    \errmessage{(Inkscape) Transparency is used (non-zero) for the text in Inkscape, but the package 'transparent.sty' is not loaded}%
    \renewcommand\transparent[1]{}%
  }%
  \providecommand\rotatebox[2]{#2}%
  \newcommand*\fsize{\dimexpr\f@size pt\relax}%
  \newcommand*\lineheight[1]{\fontsize{\fsize}{#1\fsize}\selectfont}%
  \ifx\svgwidth\undefined%
    \setlength{\unitlength}{811.53533898bp}%
    \ifx\svgscale\undefined%
      \relax%
    \else%
      \setlength{\unitlength}{\unitlength * \real{\svgscale}}%
    \fi%
  \else%
    \setlength{\unitlength}{\svgwidth}%
  \fi%
  \global\let\svgwidth\undefined%
  \global\let\svgscale\undefined%
  \makeatother%
  \begin{picture}(1,0.37448287)%
    \lineheight{1}%
    \setlength\tabcolsep{0pt}%
    \put(0,0){\includegraphics[width=\unitlength,page=1]{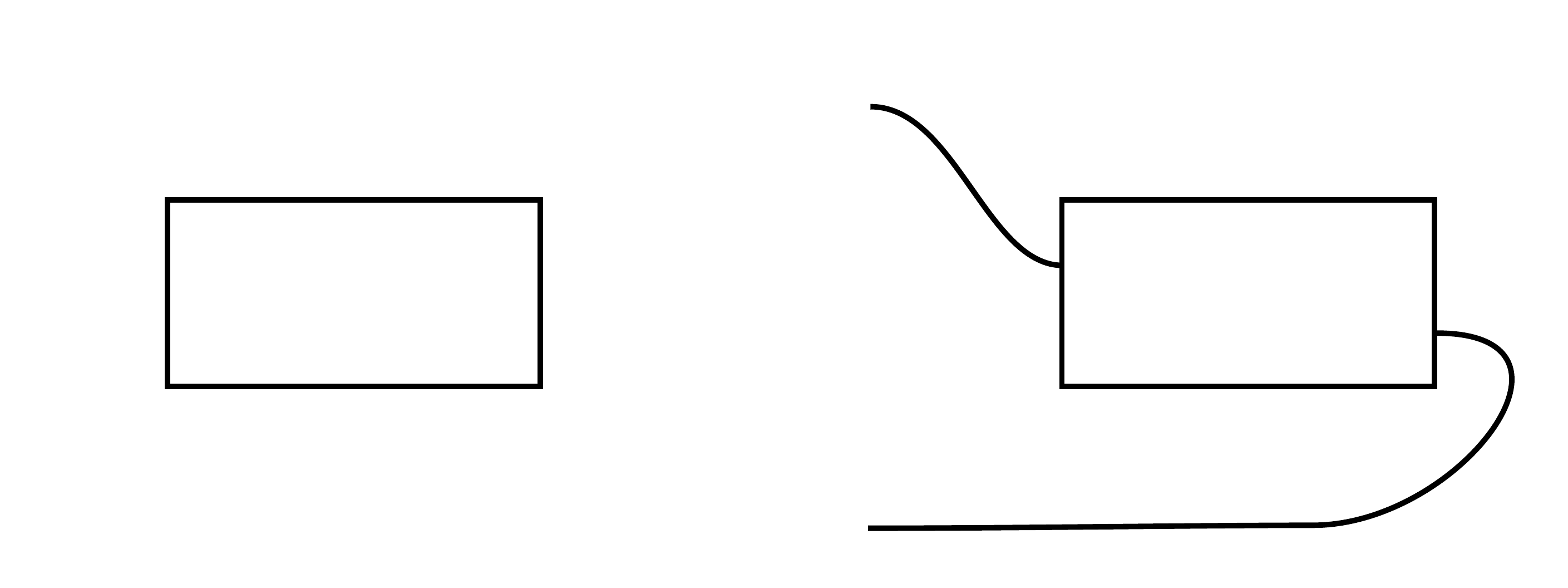}}%
    \put(0.22048816,0.15652548){\color[rgb]{0,0,0}\makebox(0,0)[t]{\lineheight{1.25}\smash{\begin{tabular}[t]{c}$B$\end{tabular}}}}%
    \put(0.79545649,0.15682331){\color[rgb]{0,0,0}\makebox(0,0)[t]{\lineheight{1.25}\smash{\begin{tabular}[t]{c}$B^{-1}$\end{tabular}}}}%
    \put(0,0){\includegraphics[width=\unitlength,page=2]{braids.pdf}}%
    \put(0.48280325,0.27899875){\color[rgb]{0,0,0}\makebox(0,0)[t]{\lineheight{1.25}\smash{\begin{tabular}[t]{c}$\Delta_1$\end{tabular}}}}%
    \put(0,0){\includegraphics[width=\unitlength,page=3]{braids.pdf}}%
    \put(0.48280325,0.0479552){\color[rgb]{0,0,0}\makebox(0,0)[t]{\lineheight{1.25}\smash{\begin{tabular}[t]{c}$\Delta_2$\end{tabular}}}}%
  \end{picture}%
\endgroup%
      \caption{Harder unknots}
      \label{F:braids}
 \end{figure}

We generalized this approach to braids with more strands in the braids, in order to hopefully increase the number of additional crossings needed.
The framework pictured in \Cref{F:braids} seems promising, where $\Delta_1$ and $\Delta_2$ are the analogues of flypes, e.g., with standard braid notation, $\Delta_1=\sigma_3\sigma_2\sigma_1\sigma_3\sigma_2\sigma_3$, see Figure~\ref{F:flype}, and $\Delta_2$ is the inverse braid.
It can be readily generalized to higher number of strands.

This construction hinges on a good choice of a four-strand pseudo-Anosov braid $B$, which can presumably not be ``flipped'' easily. Our approach is to pick one that maximizes entropy per generator, i.e., the braid $\sigma_1 \sigma_2^{-1} \sigma_3 \sigma_2^{-1}$ (see for example Thiffeault~\cite{thiffeault}).
Using the fourth power of this braid, this leads to the discovery of $D_{43}$, see \Cref{sec:hard} and \Cref{fig:harder_example}. Proving lower bounds for generalizations of this example as the number of strands goes to infinity is the object of ongoing work.


\bibliographystyle{plain}
\bibliography{biblio}

\begin{thebibliography}{10}

\bibitem{AgolCollection}
Ian Agol.
\newblock Unknots.
\newblock
  \href{https://web.archive.org/web/20120719023259/http://homepages.math.uic.edu/~agol/unknots.html}{\texttt{https://\allowbreak
  web.archive.org/\allowbreak web/\allowbreak 20120719023259/\allowbreak
  http://homepages.math.uic.edu/\allowbreak \~{}agol/\allowbreak
  unknots.html}}.

\bibitem{regina}
Benjamin~A. Burton, Ryan Budney, William Pettersson, et~al.
\newblock Regina: Software for low-dimensional topology.
\newblock \href{http://regina-normal.github.io/}{\tt http://\allowbreak
  regina-normal.\allowbreak github.\allowbreak io/}, 1999--2017.

\bibitem{CurtisCollection}
Fred Curtis.
\newblock Unknot equivalence.
\newblock
  \href{http://f2.org/maths/kt/unknoteq.html}{\texttt{http://f2.org/\allowbreak
  maths/\allowbreak kt/\allowbreak unknoteq.html}}, 2001.

\bibitem{Dynnikov06ArcPresentation}
Ivan~A. Dynnikov.
\newblock Arc-presentations of links: Monotonic simplification.
\newblock {\em Fund. Math.}, 190:29--76, 2006.

\bibitem{Freedman1994energy}
Michael~H. Freedman, Zheng-Xu He, and Zhenghan Wang.
\newblock Mobius energy of knots and unknots.
\newblock {\em Annals of Mathematics}, 139(1):1--50, 1994.

\bibitem{Goeritz34}
Lebrecht Goeritz.
\newblock Bemerkungen zur knotentheorie.
\newblock {\em Abh. Math. Sem. Univ. Hamburg}, 10(1):201--210, 1934.

\bibitem{MathOverflow_1}
Timothy Gowers~et al.
\newblock Are there any very hard unknots?
\newblock \href{https://mathoverflow.net/questions/53471/are-there-any
  -very-hard-unknots}{\texttt{https://\allowbreak mathoverflow.net/\allowbreak
  questions/\allowbreak 53471/\allowbreak are-there-any\allowbreak
  -very\allowbreak -hard-unknots}}, 2011--2017.

\bibitem{Hass01Bound}
Joel Hass and Jeffrey~C. Lagarias.
\newblock The number of {R}eidemeister moves needed for unknotting.
\newblock {\em J. Amer. Math. Soc.}, 14(2):399--428, 2001.

\bibitem{Hass08BoundBetweenDiagrams}
Joel Hass and Tahl Nowik.
\newblock Invariants of knot diagrams.
\newblock {\em Math. Ann.}, 342(1):125--137, 2008.

\bibitem{Hass10QuadraticBound}
Joel Hass and Tahl Nowik.
\newblock Unknot diagrams requiring a quadratic number of {R}eidemeister moves
  to untangle.
\newblock {\em Discrete Comput. Geom.}, 44(1):91--95, 2010.

\bibitem{Henrich14}
Allison Henrich and Louis~H. Kauffman.
\newblock Unknotting unknots.
\newblock {\em Amer. Math. Monthly}, 121(5):379--390, 2014.

\bibitem{kauffman2012hard}
Louis~H Kauffman and Sofia Lambropoulou.
\newblock Hard unknots and collapsing tangles.
\newblock {\em Introductory lectures on knot theory, Ser. Knots Everything},
  46:187--247, 2012.

\bibitem{Lackenby13PolyBoundReidemeister}
Marc Lackenby.
\newblock A polynomial upper bound on {R}eidemeister moves.
\newblock {\em Annals of Mathematics}, 182:491--564, 2015.

\bibitem{MathOverflow_2}
Daniel Moskovich~et al.
\newblock What is the state of the art for algorithmic knot simplification?
\newblock \href{https://mathoverflow.net/questions/
  144158/what-is-the-state-of-the-art-for-algorithmic
  -knot-simplification}{\texttt{https://mathoverflow.net/\allowbreak
  questions/\allowbreak 144158/\allowbreak what\allowbreak -is\allowbreak
  -the\allowbreak -state-of-the-art\allowbreak -for-algorithmic\allowbreak
  -knot-simplification}}.

\bibitem{Ochiai90}
Mitsuyuki Ochiai.
\newblock Nontrivial projections of the trivial knot.
\newblock {\em Ast\'{e}risque}, pages 7--10, 1990.
\newblock Algorithmique, topologie et g\'{e}om\'{e}trie alg\'{e}briques
  (Seville, 1987 and Toulouse, 1988).

\bibitem{petronio2016algorithmic}
Carlo Petronio and Adolfo Zanellati.
\newblock Algorithmic simplification of knot diagrams: New moves and
  experiments.
\newblock {\em Journal of Knot Theory and Its Ramifications},
  25(10):1650059:1--30, 2016.

\bibitem{Reidemeister27}
Kurt Reidemeister.
\newblock Elementare {B}egr\"{u}ndung der {K}notentheorie.
\newblock {\em Abh. Math. Sem. Univ. Hamburg}, 5(1):24--32, 1927.

\bibitem{sagemath}
{The Sage Developers}.
\newblock {\em {S}ageMath, the {S}age {M}athematics {S}oftware {S}ystem
  ({V}ersion 9.0)}, 2020.
\newblock \href{https://www.sagemath.org}{\tt https://www.sagemath.org}.

\bibitem{thiffeault}
Jean-Luc Thiffeault and Matthew~D Finn.
\newblock Topology, braids and mixing in fluids.
\newblock {\em Philosophical Transactions of the Royal Society A: Mathematical,
  Physical and Engineering Sciences}, 364(1849):3251--3266, 2006.

\bibitem{Tuzun16}
Robert~E. Tuzun and Adam~S. Sikora.
\newblock Verification of the {J}ones unknot conjecture up to 22 crossings.
\newblock {\em J. Knot Theory Ramifications}, 27(3):1840009, 18, 2018.

\end{thebibliography}


\addresseshere

\newpage

\appendix

\section{Gauss codes of hard unknot diagrams}
\label{appendix}

Here we give Gauss codes and planar drawings\footnote{all drawings of the section are generated with {SAGE}~\cite{sagemath}} for all unknot diagrams featuring in this note. We start by repeating the table from \Cref{sec:errors}, now including $D_{28}$, $D_{43}$ and $PZ_{78}$.

\bigskip

\begin{center}
\begin{tabular}{llll}
  \toprule
  Name & $\crs$ & $\extra(D)$ & References \\
  \midrule
  $D_{28}$ 		&$28$   &$3$ 	&	\S\ref{sec:hard} \\
  $D_{43}$ 		&$43$   &$\geq 3$&	\S\ref{sec:hard} \\
  $PZ_{78}$ 	&$78$ &$\geq 3$	&   \S\ref{sec:hard} \\
  \midrule
  H	 		&$9$   &$1$ 	&	\cite[Figure~4]{Henrich14} \\
  J	 		&$9$   &$1$ 	&	\cite[Figure~4]{Henrich14} \\
  Culprit 		&$10$  &$1$	&	\cite{kauffman2012hard} \\
  Monster 		&$10$  &$0$	&	\cite{petronio2016algorithmic, CurtisCollection}\\
  Goeritz 		&$11$  &$1$	&	\cite{Goeritz34}, cf. \S\ref{ssec:goeritz} \\
  Thistlethwaite 	&$15$  &$0$	& \cite[Figure~9]{petronio2016algorithmic}, \cite{CurtisCollection} \\
  Ochiai I 		&$16$  &$0$	&	\cite[Figure 1]{Ochiai90} \\
  Tuzun-Sikora		&$21$  &$0$	&	\cite[Figure 8]{Tuzun16} \\
  Freedman-He-Wang 		&$32$  &$0$	&	\cite[Figure 6.1]{Freedman1994energy}, cf. \S\ref{ssec:freedman} \\
  ``Fake'' Freedman-He-Wang 		&$32$  &$0$	&	\S\ref{ssec:freedman} \\
  Ochiai II 		&$45$  &$0$	&	\cite[Figure~2]{Ochiai90} \\
  Ochiai II (reduced) 	&$35$  &$\geq 2$&	\S\ref{sec:errors} \\
  Ochiai III 	&$67 $  &$ 0$&		\cite[Figure~3]{Ochiai90} \\
  Ochiai IV (Suzuki) 		&$55$  &$\geq 2$&	\cite[Figure~4]{Ochiai90} \\
  Haken 		&$141$ &$0$	&	 \cite{petronio2016algorithmic} \\
  $PZ_{31}$ 	&$31$ &$0$	&	\cite[Figure~12]{petronio2016algorithmic} \\
  $PZ_{120}$	&$120$ &$\geq 2$	&	\cite[Figure~27]{petronio2016algorithmic} \\
  $PZ_{138}$	&$138$ &$2$	&	\cite[Figure~14]{petronio2016algorithmic} \\
  \bottomrule
\end{tabular}
\end{center}

\bigskip

\begin{figure}[h]
  \centerline{\includegraphics[height=6cm]{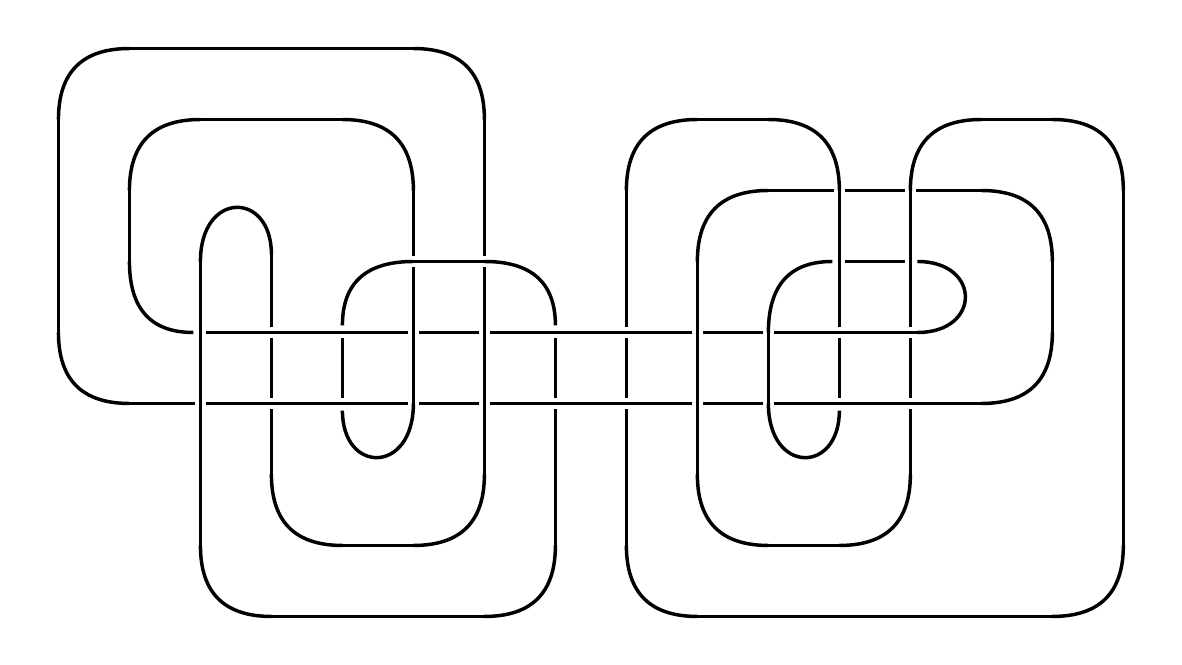}}
  \caption{A $28$-crossing diagram $D_{28}$ of the unknot requiring three extra crossings. \\ Gauss code: 1 -4 -3 6 5 -2 -7 8 4 -5 -9 10 2 -1 -11 7 12 -13 -6 3 14 -12 -10 9 13 -14 -8 11 15 -18 -17 20 19 -16 -21 17 22 -23 -20 21 24 -15 -25 26 16 -19 -27 28 18 -24 -26 27 23 -22 -28 25 \label{fig:hard_example}}
\end{figure}

\begin{figure}
  \centerline{\includegraphics[height=9cm]{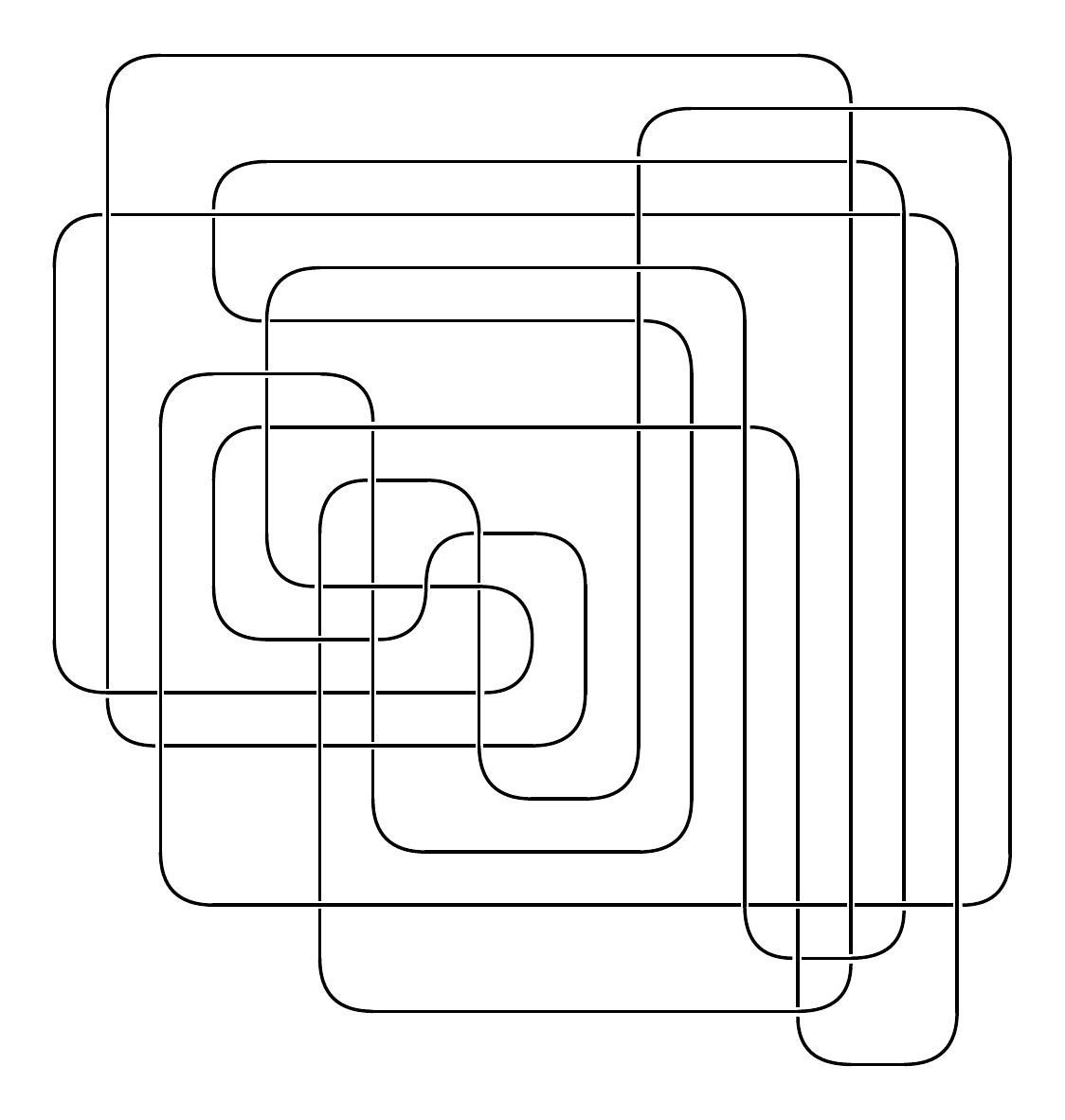}}
  \caption{A $43$-crossing diagram $D_{43}$ of the unknot requiring at least three extra crossings. \\ Gauss code: \texttt{-1 2 -3 -4 -5 6 -7 8 -9 10 -11 12 -13 3 14 -15 16 -17 18 -19 20 1 -21 22 -23 24 -25 7 -26 27 -10 28 -29 13 4 30 31 -32 17 -33 34 35 -2 -36 -22 -37 -38 25 -8 39 -27 11 -40 29 41 23 37 -31 42 -16 33 -43 19 -20 -35 -14 -30 5 -6 38 -24 9 -39 26 -12 40 -28 -41 36 21 32 -42 15 -34 43 -18} \label{fig:harder_example}}
\end{figure}

\begin{figure}[h]
  \centerline{\includegraphics[width=\textwidth]{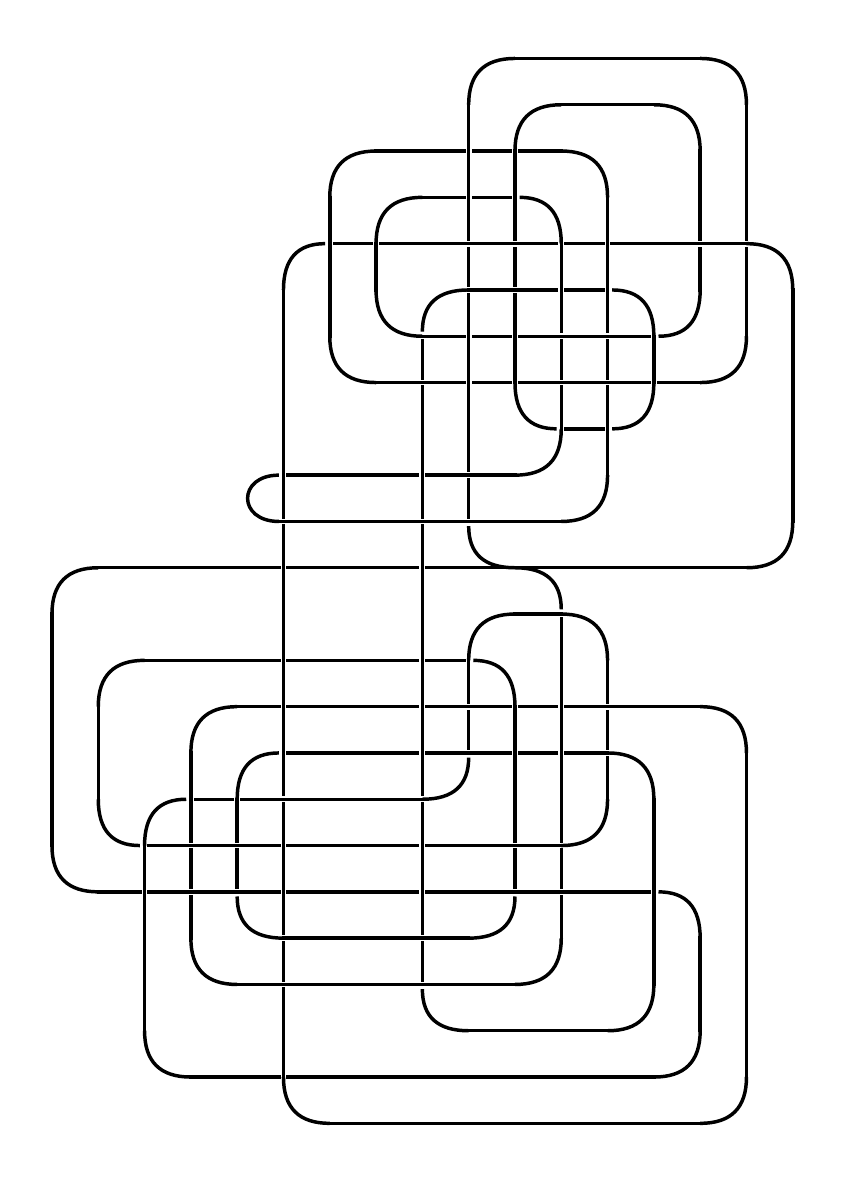}}
  \caption{A $78$-crossing diagram $PZ_{78}$ of the unknot obtained from reducing $138$-crossing diagram $PZ_{138}$ from~\cite[Figure~14]{petronio2016algorithmic} via a $140$-crossing diagram. \\ Gauss code: \texttt{1 2 -3 -4 5 6 -7 -8 9 10 -2 -11 12 13 -14 -5 15 16 17 18 -19 -20 21 22 -6 -23 24 25 -26 -1 27 28 -29 -30 23 14 -31 -32 11 26
 -33 -34 30 7 -22 -15 4 31 -13 -24 34 36 -28 -9 20 -17 37 38 35 39 -40 41 42 -43 -44 45 46 -47 -48 49 -39 -50 51 -52 -53 54 43 -55 -56 48 57 -58 -38 -61 -59 60 52 -62 -41 56 63 -46 -64 65 58 -66 -49 40 67 -51 -68 59 69 -70 -45 71 55 -42 -72 53 74 -69 -73 -37 -65 75 47 -63 -71 44 76 -74 -60 68 -77 -35 66 -57 -75 64 70 -76 -54 72 62 -67 50 77 61
 73 -18 19 -10 -27 78 33 -25 -12 32 3 -16 -21 8 29 -36 -78} \label{fig:PZ_example}}
\end{figure}

\newpage

\begin{figure}[h]
 \centerline{\includegraphics[height=7cm]{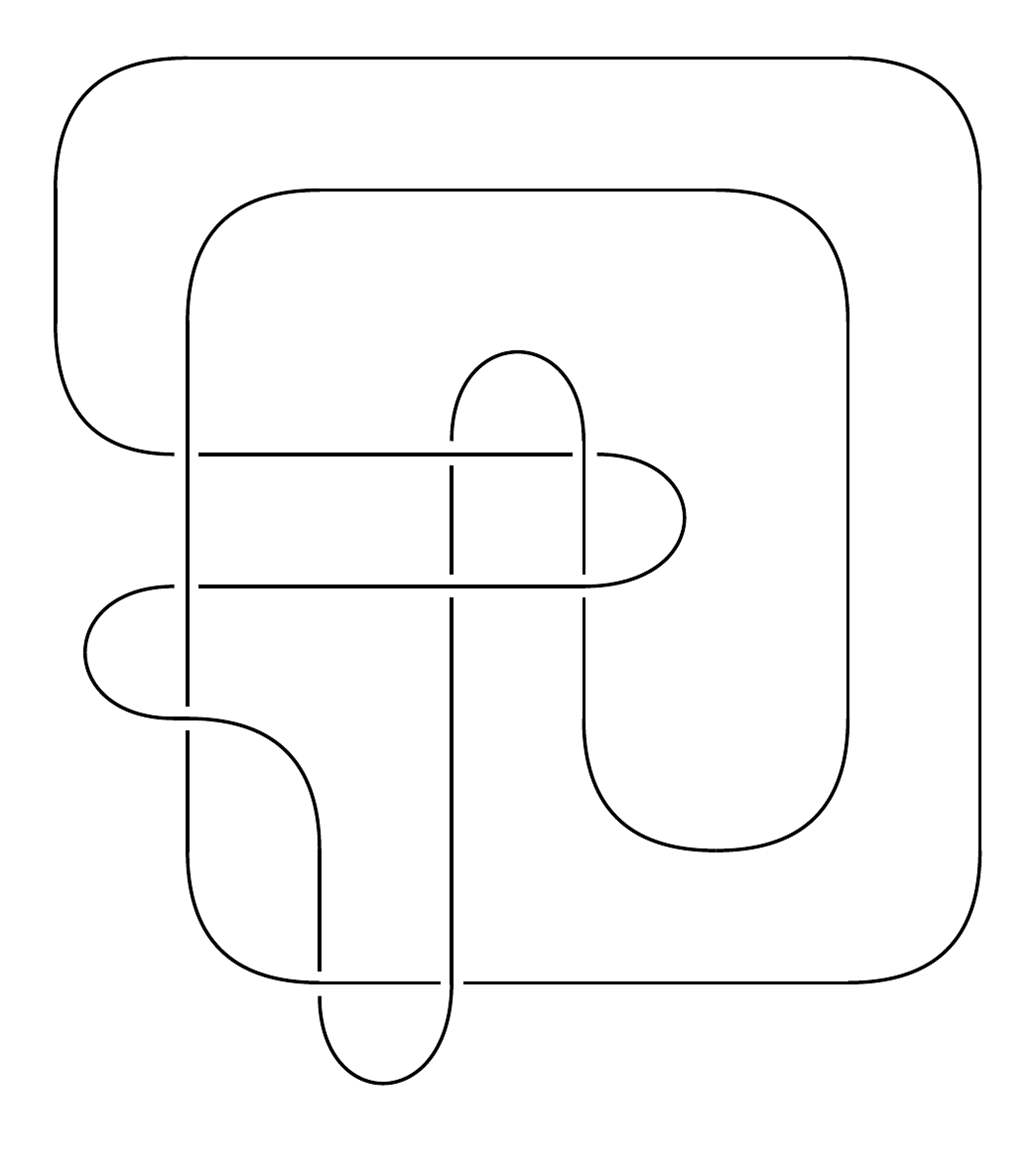}}
 \caption{{\em H (\cite[Fig.~4]{Henrich14})}: \texttt{-1 2 -3 4 -5 -6 7 -8 9 1 -2 3 -4 -9 6 -7 8 5}
}
\end{figure}

\begin{figure}[h]
 \centerline{\includegraphics[height=7cm]{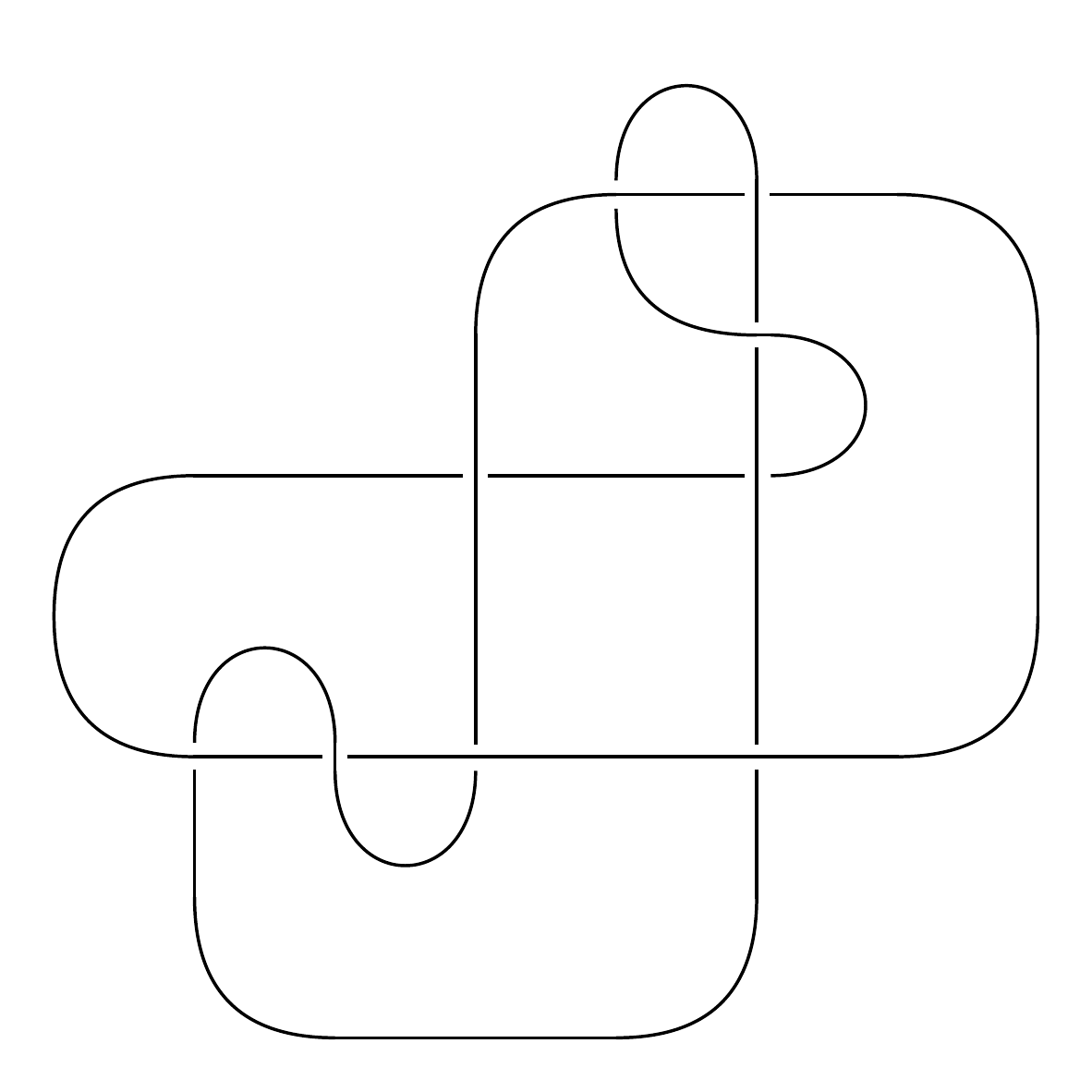}}
 \caption{{\em J (\cite[Fig.~4]{Henrich14})}: \texttt{-1 2 -3 4 5 -6 7 1 -4 3 -2 -5 8 -9 6 -7 9 -8}
}
\end{figure}

\begin{figure}[h]
 \centerline{\includegraphics[height=7cm]{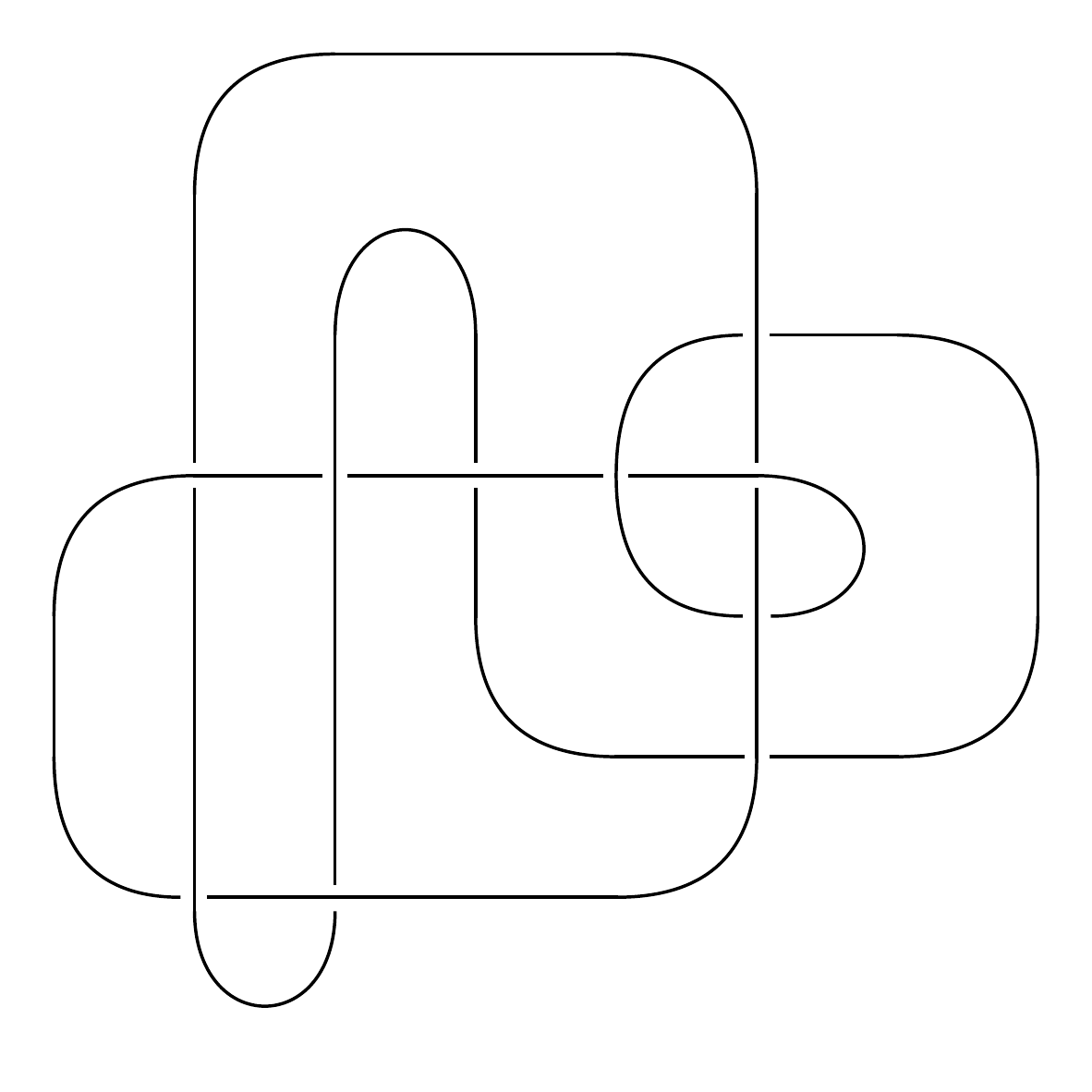}}
 \caption{{\em Culprit (\cite[Figure~2, Figure~15]{Henrich14})}: \texttt{-1 2 -3 4 -5 6 7 8 -9 10 -4 5 -6 3 -2 -7 -10 1 -8 9}
}
\end{figure}

\begin{figure}[h]
 \centerline{\includegraphics[height=7cm]{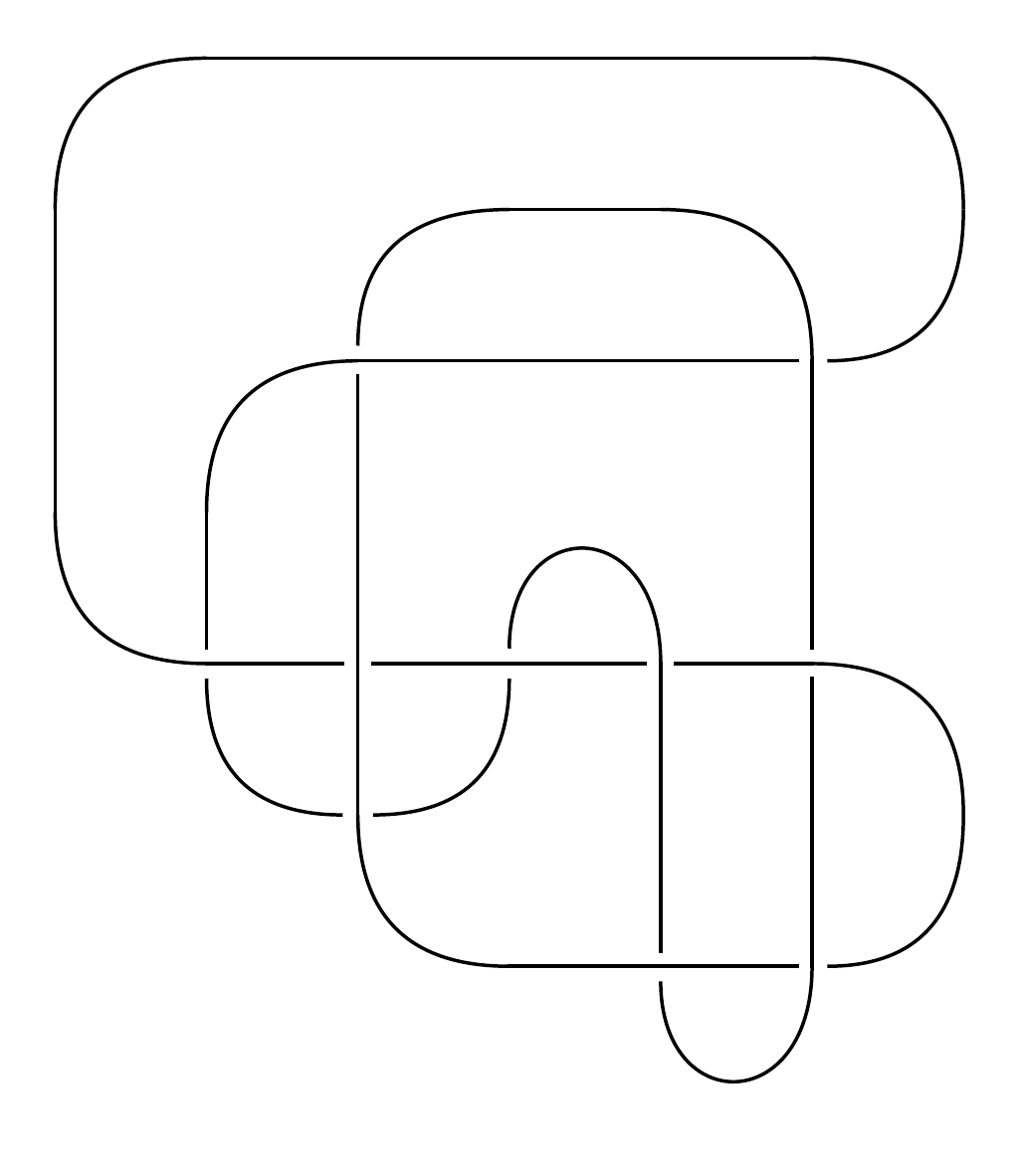}}
 \caption{{\em Monster}: \texttt{ 1 -2 3 4 5 -6 7 -8 9 -3 10 -1 2 -10 -4 -9 8 -5 6 -7}}
\end{figure}

\begin{figure}[h]
 \centerline{\includegraphics[height=6cm]{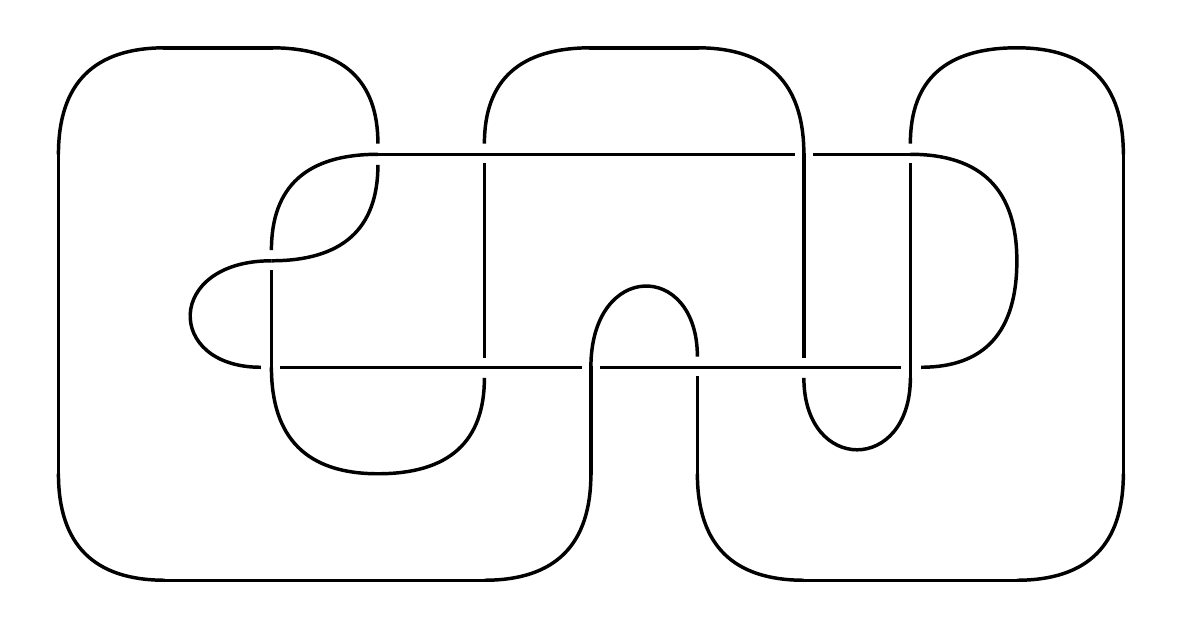}}
 \caption{{\em Goeritz (\cite{Goeritz34})}: \texttt{ 1 -2 3 -4 -5 6 -7 8 -9 -10 11 -1 2 -3 4 -11 10 7 -8 9 -6 5}}
\end{figure}

\begin{figure}[h]
 \centerline{\includegraphics[height=7cm]{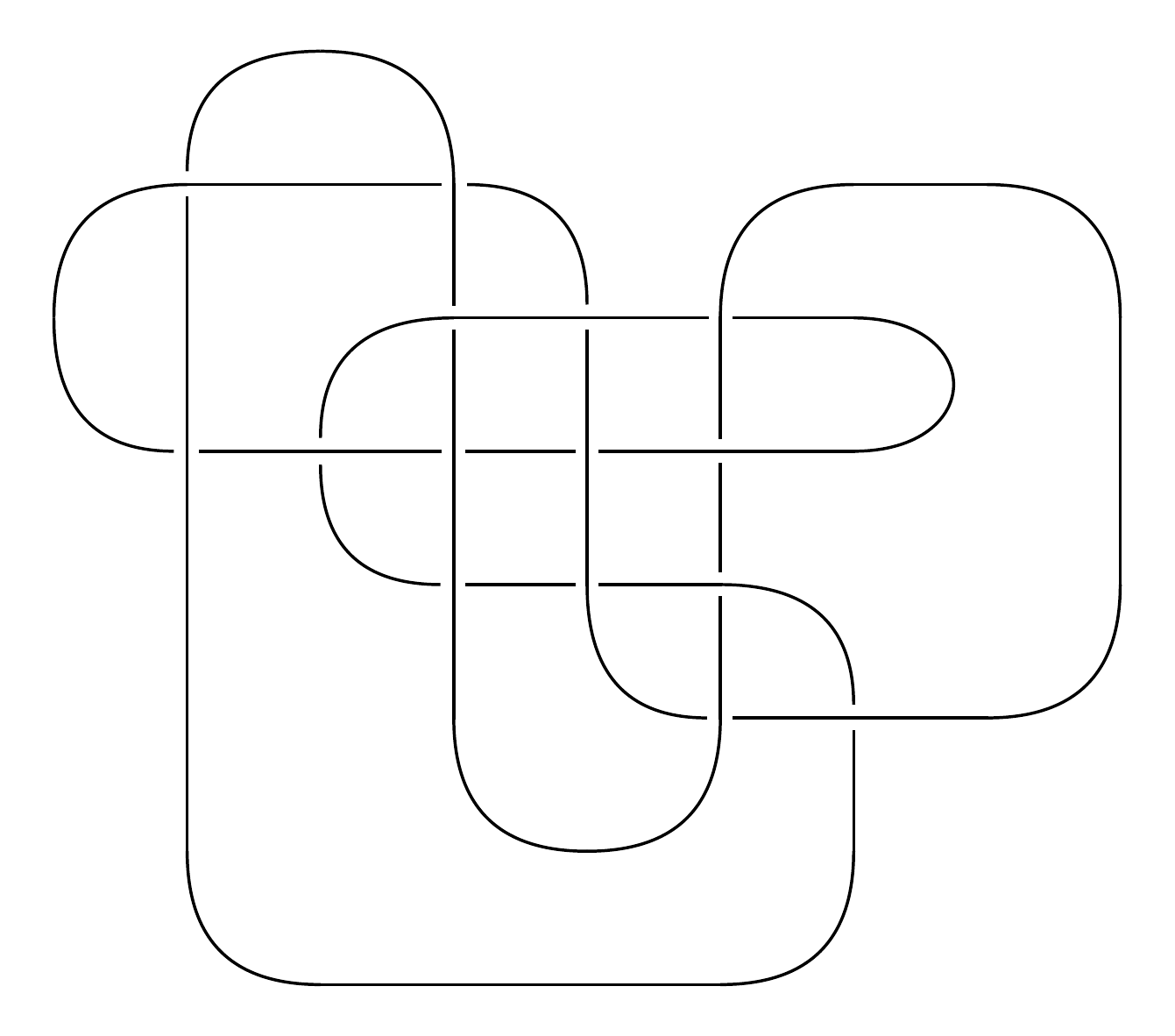}}
 \caption{{\em Thistlethwaite}: \texttt{ 1 2 -3 4 5 -6 -7 8 -9 10 -11 -5 12 -1 6 13 -10 -14 -4 15 -2 9 -8 7 -13 11 14 3 -15 -12}
}
\end{figure}

\begin{figure}[h]
 \centerline{\includegraphics[height=7cm]{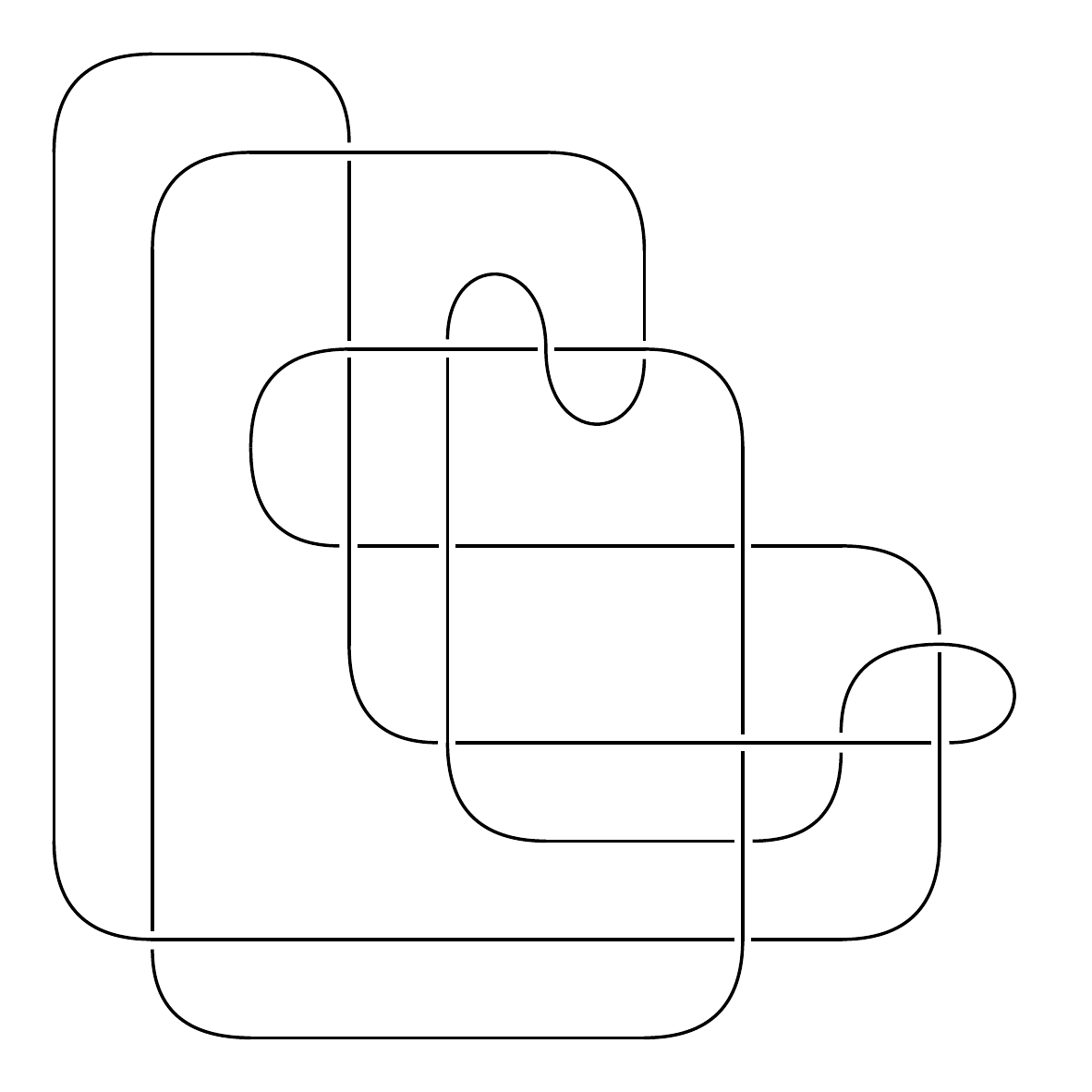}}
 \caption{{\em Ochiai I (\cite[Figure~1]{Ochiai90})}: \texttt{-1 2 3 -4 -5 -6 -7 8 -9 10 -11 -3 4 -12 13 14 -8 7 -14 -15 12 5 -2 1 -16 11 -10 9 15 -13 6 16}
}
\end{figure}

\begin{figure}[h]
 \centerline{\includegraphics[height=7cm]{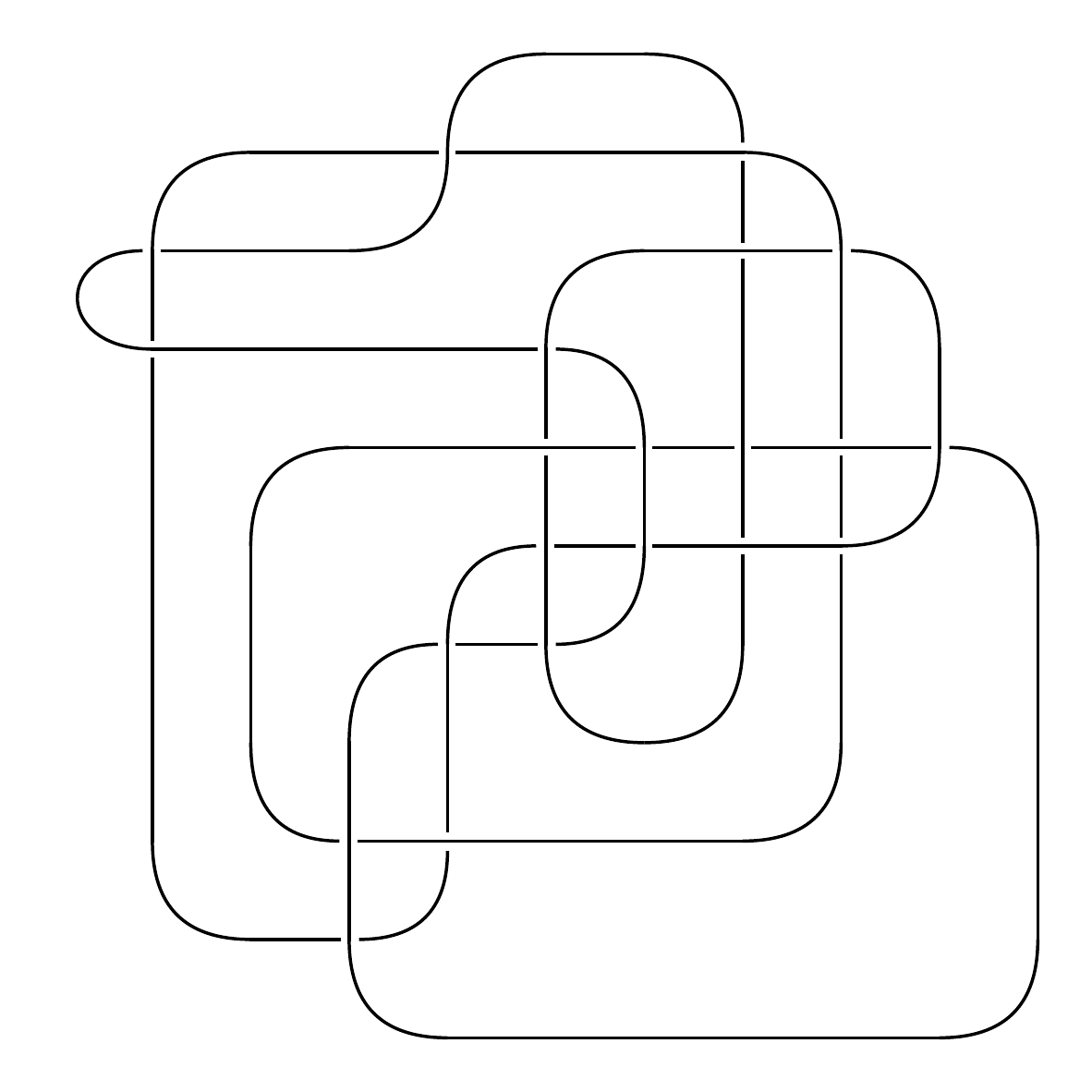}}
 \caption{{\em Tuzun Sikora example (\cite[Figure 8]{Tuzun16})}: \texttt{1 6 -7 10 11 12 -14 -15 16 -17 -18 -2 3 -4 5 7 -9 -11 17 -21 19 -20 -8 9 -10 18 21 -16 -13 14 20 -1 2 -3 4 -5 -6 8 -12 13 15 -19}
}
\end{figure}

\begin{figure}[h]
 \centerline{\includegraphics[height=7cm]{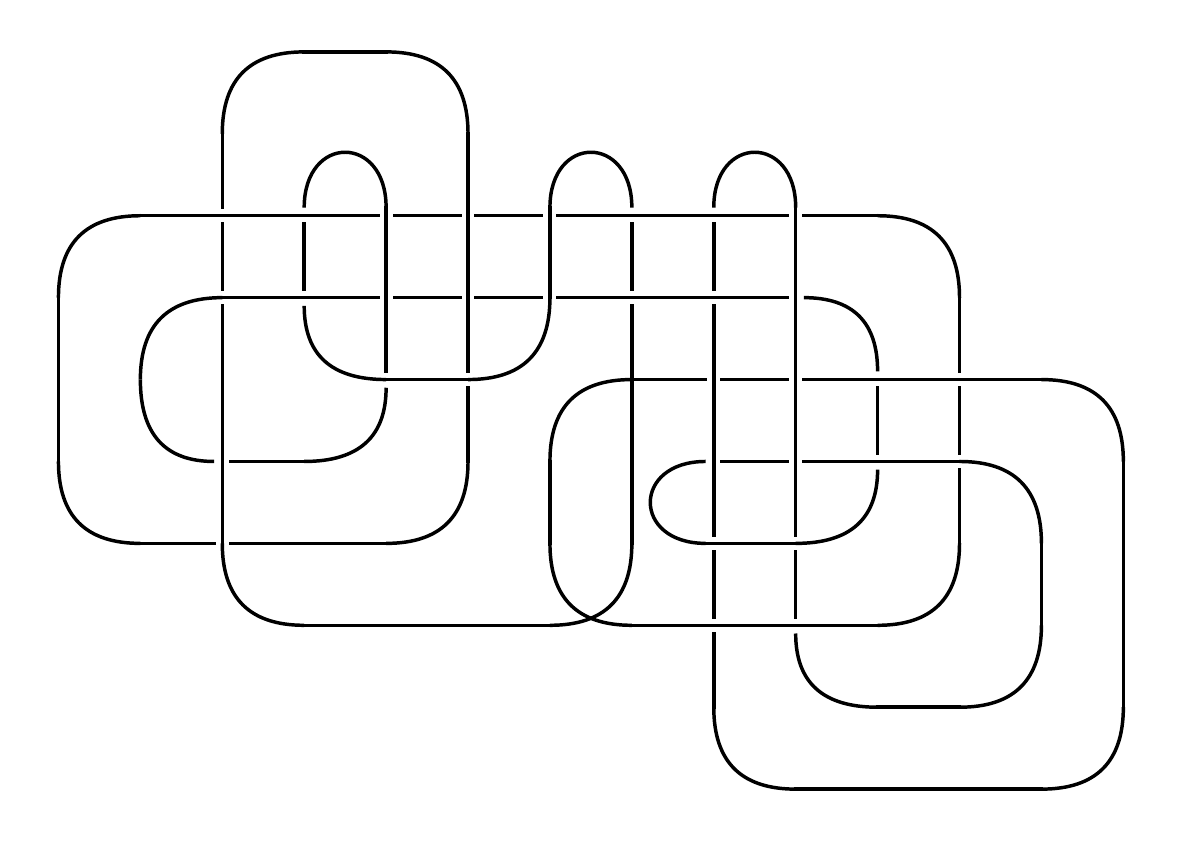}}
 \caption{{\em ``Fake'' Freedman-He-Wang}: \texttt{ 1 -2 -3 -4 5 6 -7 -8 4 9 -10 -5 8 11 2 12 -13 -1 14 7 -6 -15 16 3 -11 -14 15 10 -9 -16 -12 13 17 -18 -19 -20 21 22 -23 -24 25 19 -22 -26 27 23 -28 -17 18 29 24 30 -31 -21 20 32 -30 -27 26 31 -32 -25 -29 28} \label{fig:fakeFHWD}}
\end{figure}

\begin{figure}[h]
 \centerline{\includegraphics[height=7cm]{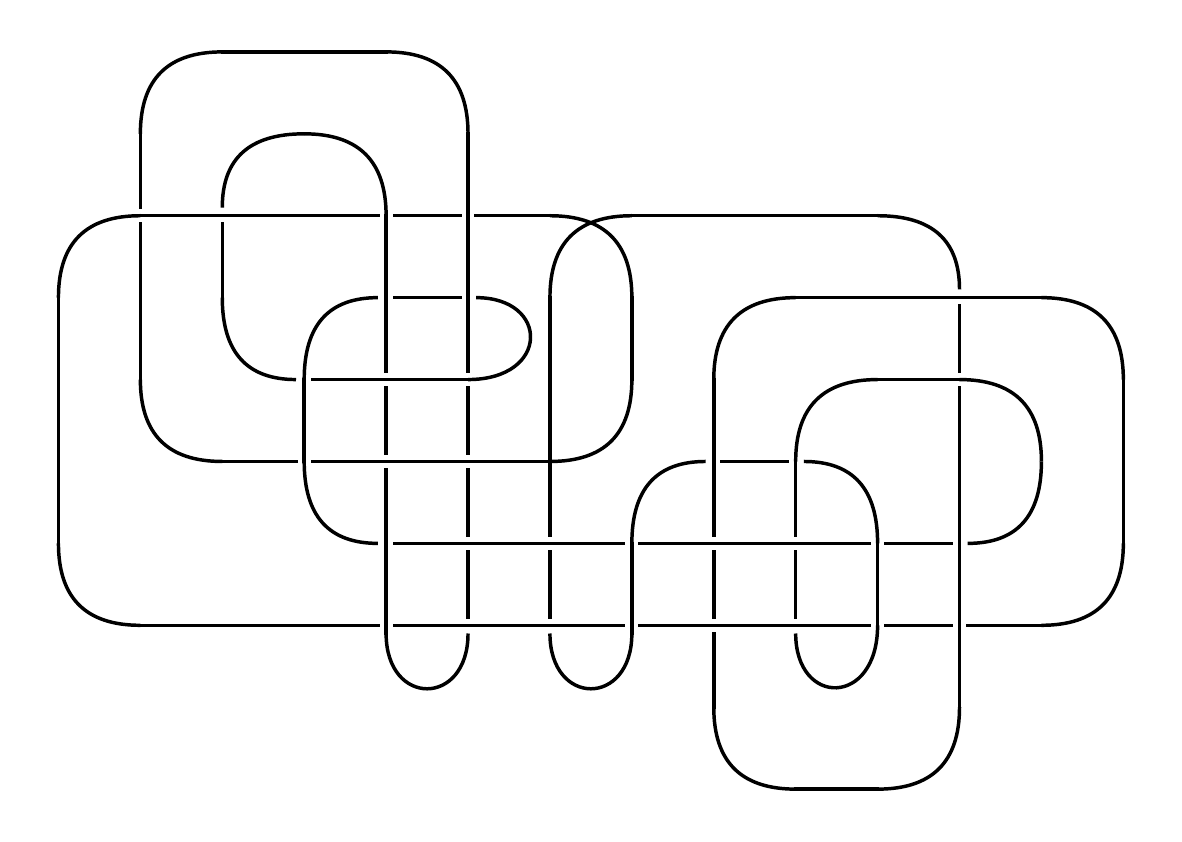}}
 \caption{{\em Freedman-He-Wang (\cite[Figure~6.1]{Freedman1994energy})}: \texttt{ 3 4 -5 -8 10 12 -13 -15 8 7 -9 -10 15 16 -4 1 17 -20 32 31 -27 -25 22 24 -31 -29 28 27 -24 -23 20 19 -18 -17 -21 -22 25 26 -30 -32 23 21 -26 -28 29 30 -19 18 2 -3 14 13 -12 -11 6 5 -16 -14 11 9 -7 -6 -1 -2} \label{F:hard_unknot}}
\end{figure}

\begin{figure}[h]
 \centerline{\includegraphics[height=8cm]{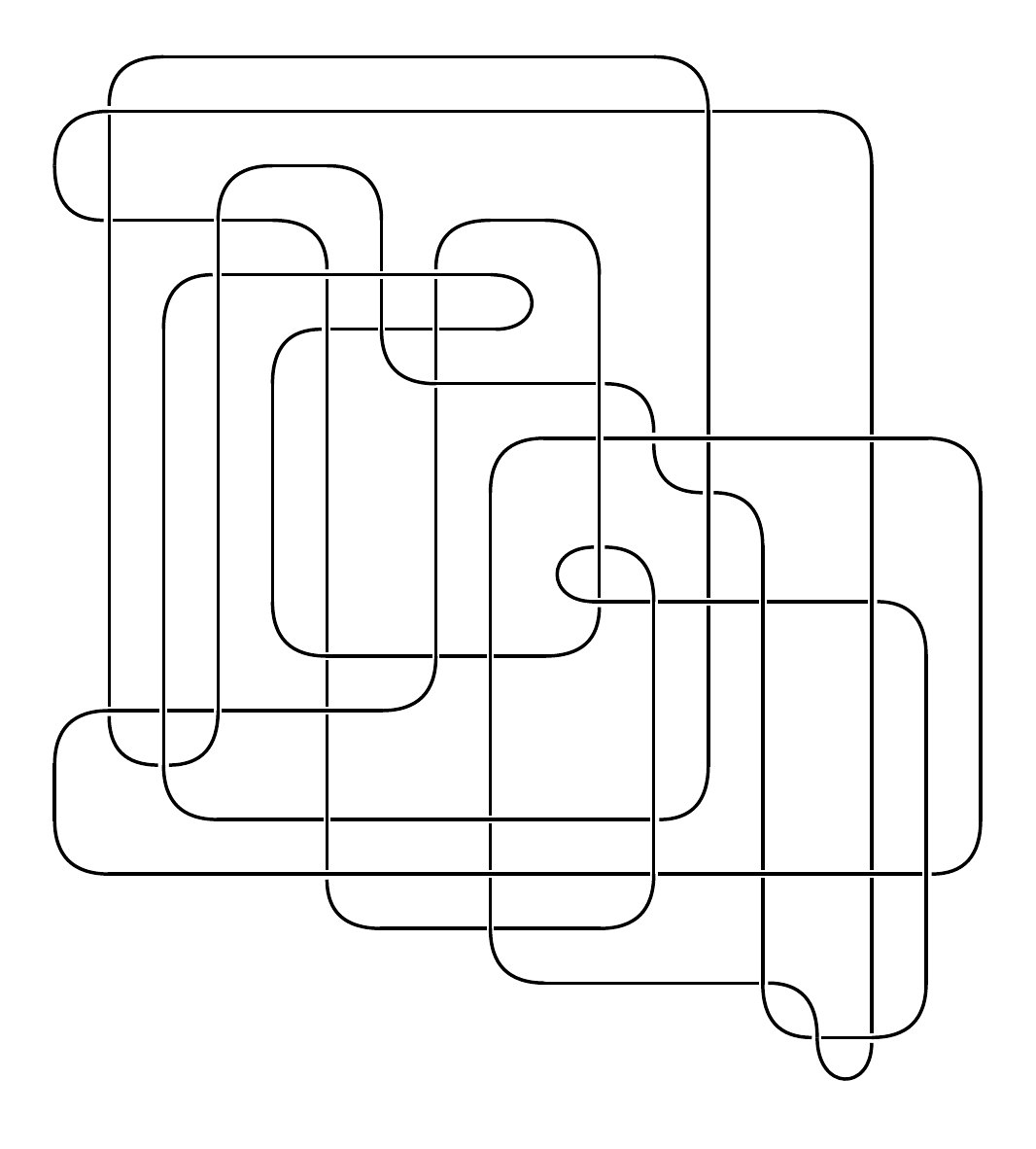}}
 \caption{{\em Ochiai II (\cite[Figure~2]{Ochiai90})}: \texttt{-1 2 -3 -4 5 -6 -7 8 -9 10 -11 -12 13 -14 15 16 4 -17 6 18 -19 11 20 21 -22 -23 24 25 -26 27 -28 29 14 30 -31 -13 -25 32 -33 -34 -27 28 35 36 -37 22 38 -39 12 -40 -41 42 -16 3 -2 1 -21 -38 23 37 45 -35 34 26 -29 -15 -42 43 -18 9 -8 7 17 -5 -43 41 -44 31 -30 44 40 19 -10 -20 39 -24 -32 33 -36 -45}
}
\end{figure}

\begin{figure}[h]
 \centerline{\includegraphics[height=6cm]{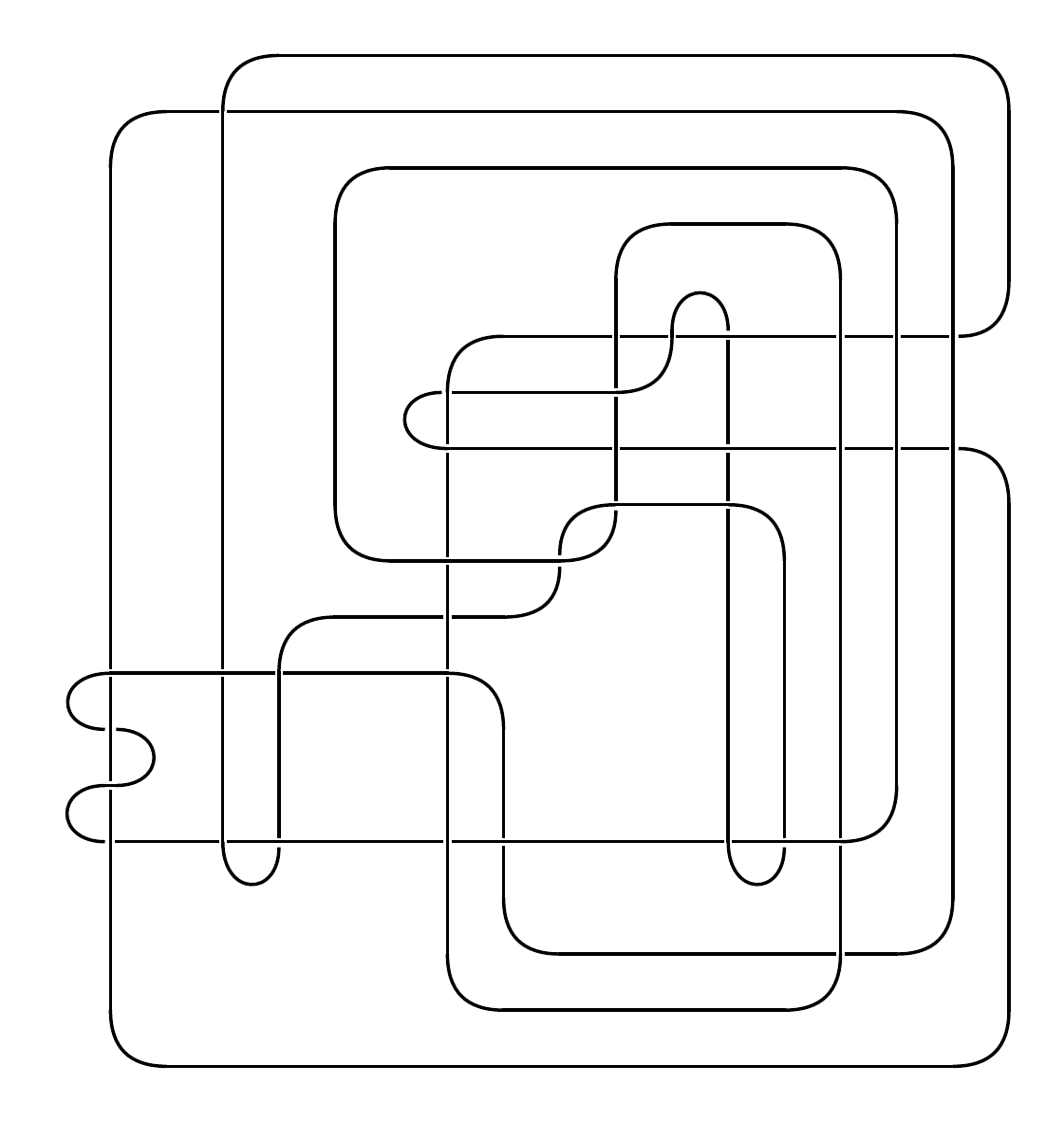}}
 \caption{{\em Ochiai II (reduced)}: \texttt{ 1 -2 -21 -3 4 -6 -5 7 -8 -10 9 -11 12 15 -13 14 3 5 -16 17 -18 19 10 21 22 13 23 -24 20 -12 25 -26 -27 28 -30 31 29 -32 11 -20 -15 33 2 -1 -31 30 -28 27 -33 -22 -14 34 -17 8 -7 16 6 -4 -34 -35 24 -23 35 18 -19 -9 32 -25 26 -29}}
\end{figure}

\begin{figure}[h]
 \centerline{\includegraphics[height=6.5cm]{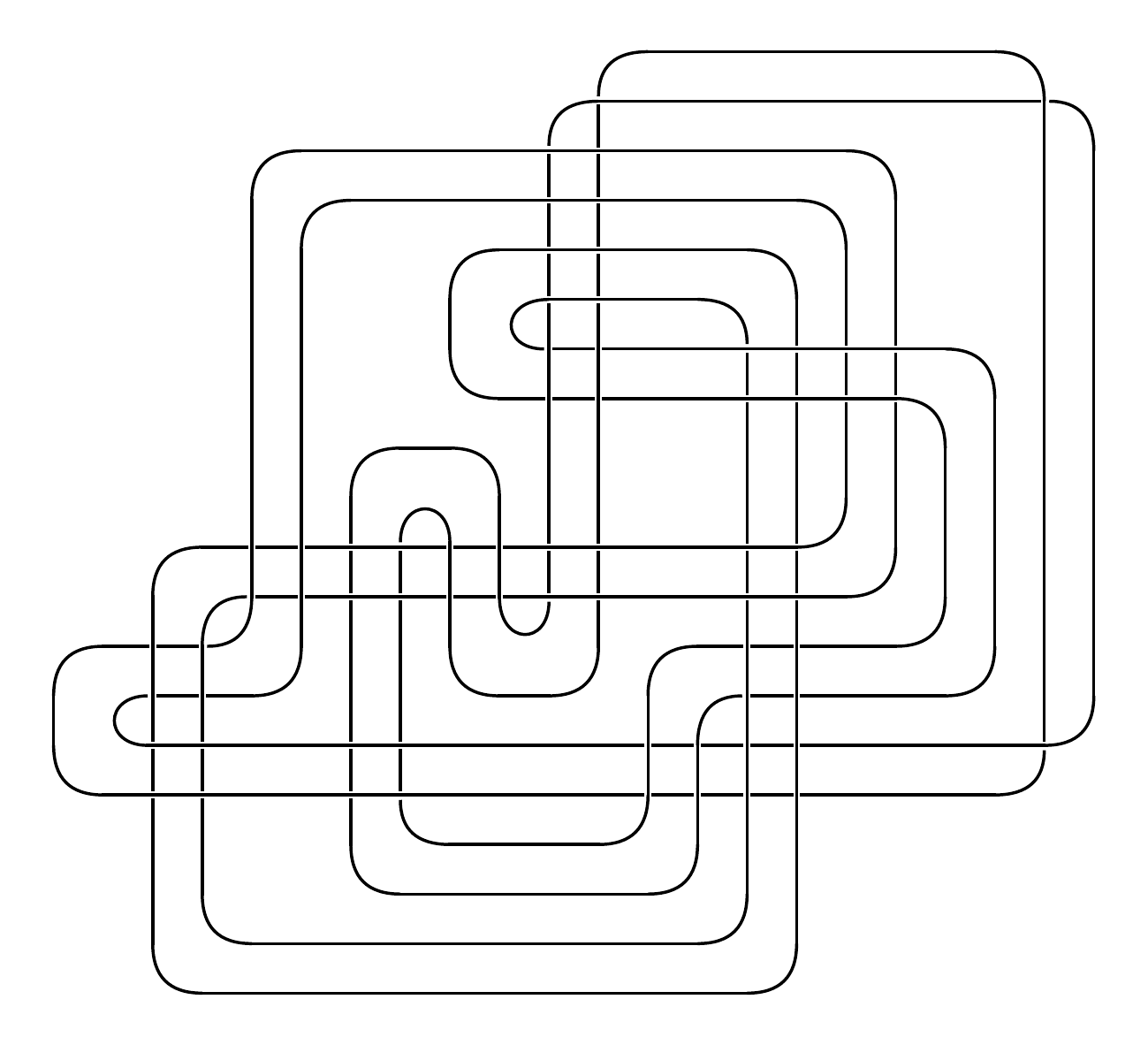}}
 \caption{{\em Ochiai III (\cite[Figure~3]{Ochiai90})}: \texttt{1 3 31 30 -45 -44 65 64 62 63 -56 -57 54 55 -3 -4 7 8 -10 -9 24 23 49 48 -67 -65 -42 -43 33 32 -38 -39 40 42 44 46 -48 -50 17 18 -13 -15 -52 -54 57 58 -61 -62 39 37 -34 -33 -30 -29 -27 26 -25 -24 -22 -19 -18 13 14 11 9 -7 -5 2 4 28 29 -47 -46 67 66 61 60 -59 -58 52 53 -1 -2 5 6
-12 -11 22 21 51 50 -66 -64 -40 -41 34 35 -36 -37 41 43 45 47 -49 -51 20 19 -14 -16 -53 -55 56 59 -60 -63 38 36 -35 -32 -31 -28 27 -26 25 -23 -21 -20 -17 15 16 12 10 -8 -6}}
\end{figure}

\begin{figure}[h]
 \centerline{\includegraphics[height=7.5cm]{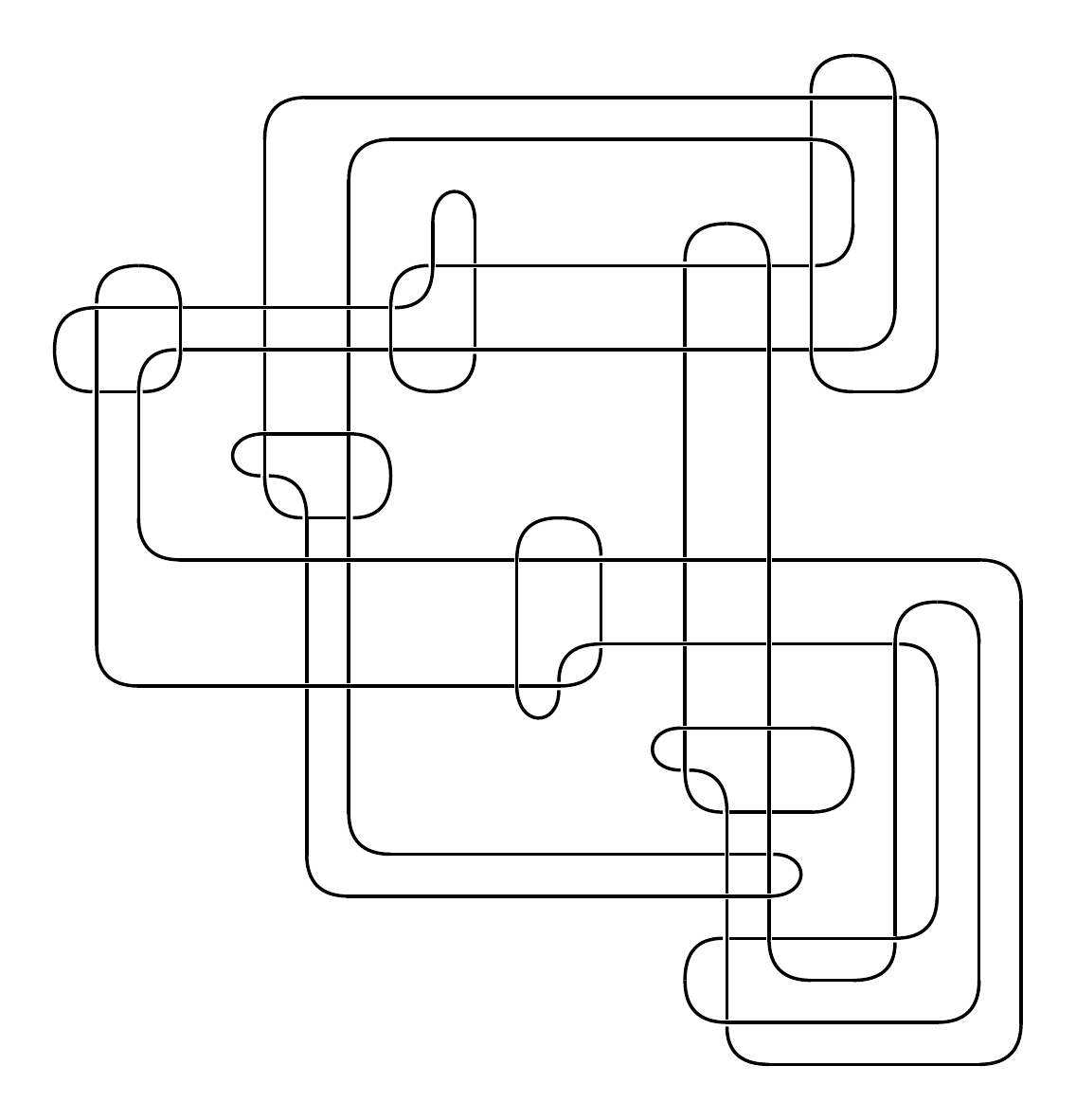}}
 \caption{{\em Ochiai IV (Suzuki,~\cite[Figure~4]{Ochiai90})}: \texttt{ 1 -2 3 6 7 -8 9 -10 -11 12 8 -9 10 13 -14 -15 16 -17 -19 20 15 -16 17 -21 22 23 -24 14 25 26 27 -28 -29 -30 -13 -31 32 -33 -23 24 31 -32 33 34 -35 19 -20 -25 30 11 -12 36 37 4 -5 38 39 -40 41 29 -26 -42 43 -44 -45 46 42 -43 44 -47 -38 -48 49 -50 -51 52 48 -49 50 -37 -6 35 21 -22 -34 -7 -36 51 -52 -39 -53 45 -46 -27 28 54 -55 18 40 -41 -54 55 -18 53 47 -1 2 -3 -4 5}}
\end{figure}

\begin{figure}[h]
 \centerline{\includegraphics[width=\textwidth]{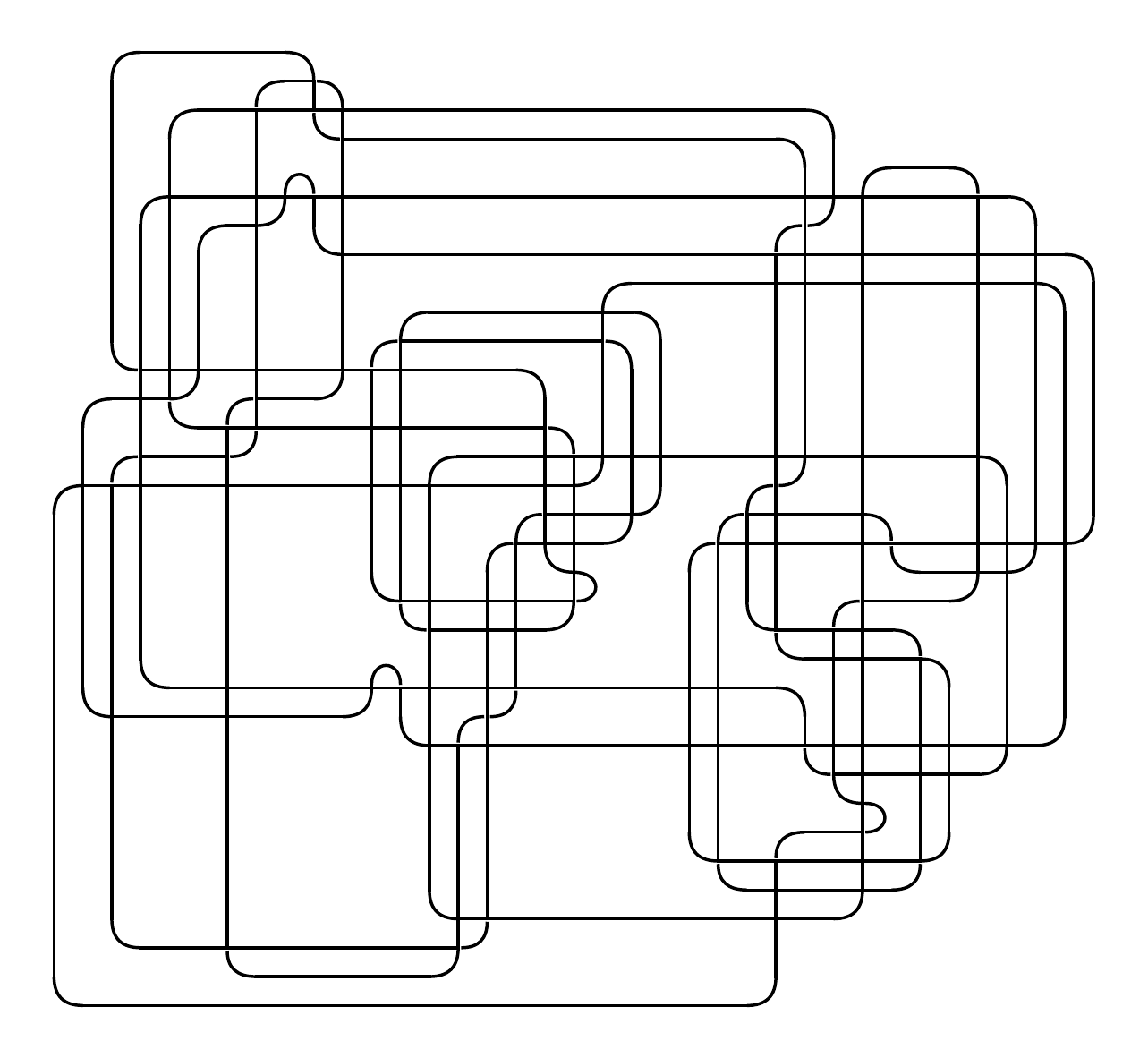}}
 \caption{{\em Haken}: \texttt{-1 -2 3 -4 5 6 7 8 -9 -10 11 12 13 -14 15 -16 17 18 -19 20 -21 -7 22 23 24 -25 -26 -27 28 -29 30 31 -32 -33 -34 35 36 -37 -38 -39 -40 -41 42 43 -44 45 46 47 -48 -49 -50 38 -51 52 -53 -54 -55 -42 -56 57 58 -46 -59 -60 -61 -62 63 -64 -65 -66 -67 68 -35 -69 70 -71 -72 73 29 -74 -31 -75 -76 77 69 -36 -68 -78 -79 -80 64 -63 -81 -82 60 -83 -47 84 -57 -85 -43 55 -86 -87 -52 -88 37 -89 49 -90 -84 -58 -45 44 85 56 -91 -92 39 51 88 -70 -77 -93 33 -94 74 -30 -73 95 27 -96 -97 98 -23 99 -8 21 100 -101 -18 102 -103 -15 -104 105 -12 106 -107 9 -99 -22 108 -5 -109 110 2 -111 -112 113 -114 -28 -95 72 115 -116 75 32 94 117 78 -118 66 -119 120 121 83 59 10 107 -122 -123 124 103 16 125 -110 -3 126 127 -98 -24 54 86 -128 40 92 129 90 48 -121 -130 62 81 -131 104 14 132 -133 -124 -134 19 101 135 -6 -108 -127 136 137 114 -113 -138 139 -126 4 109 -135 -100 -20 140 123 133 -132 -13 -105 131 141 79 118 67 89 50 -129 91 41 128 87 53 25 97 -136 -139 111 1 -125 -17 -102 134 -140 122 -106 -11 82 61 130 -120 119 65 80 -141 -117 34 93 76 116 -115 71 26 96 -137 138 112}}
\end{figure}

\begin{figure}[h]
  \centerline{\includegraphics[width=\textwidth]{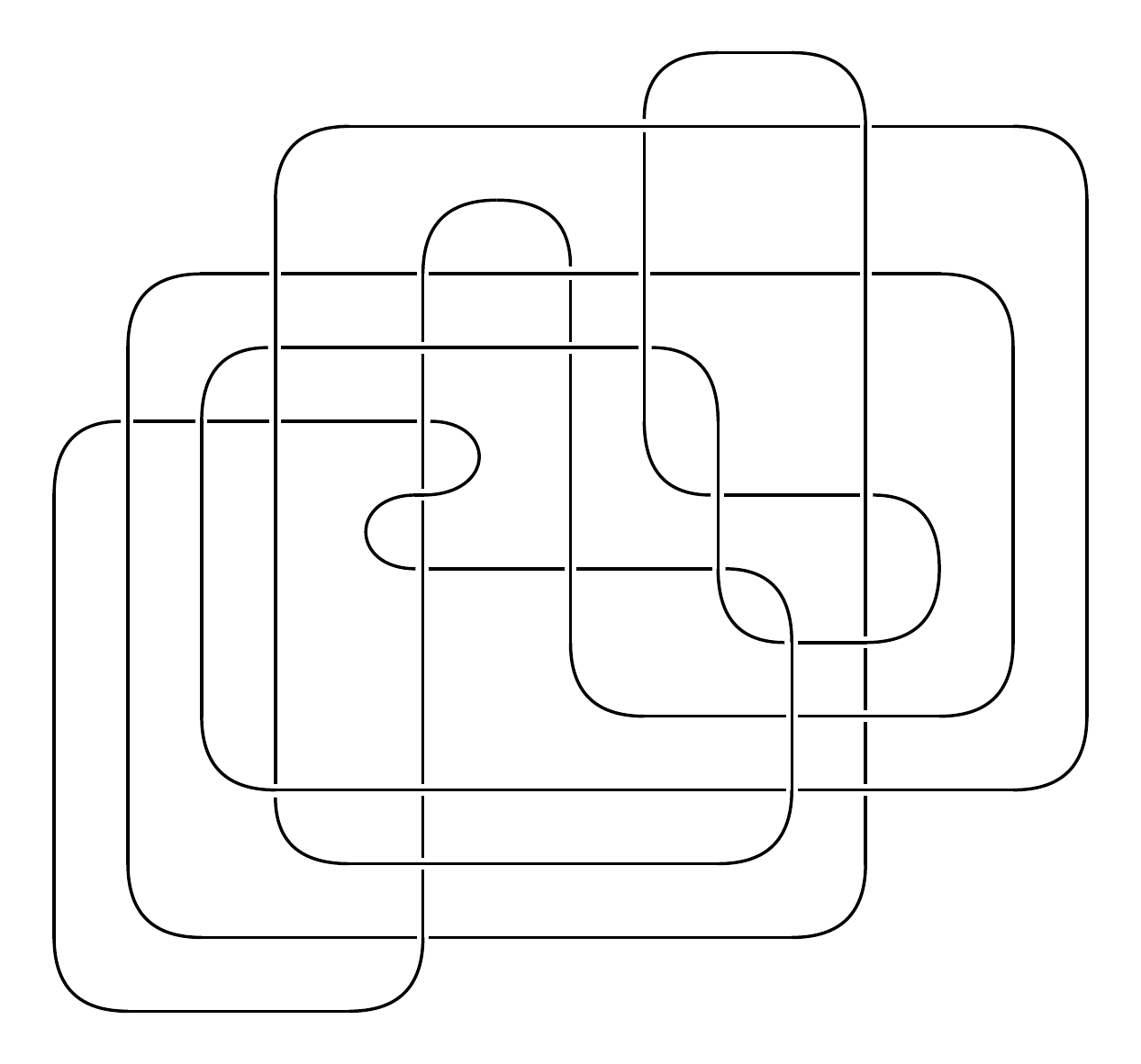}}
  \caption{$PZ_{31}$ {\em (\cite[Figure~12]{petronio2016algorithmic})}: \texttt{2 -3 4 -19 16 -21 -22 1 -20 26 -27 -29 21 -16 -17 7 -8 -13 -26 -25 28 27 14 -15 29 30 -31 -28 25 -24 23 31 -30 22 19 -18 6 9 11 -12 13 -14 15 17 18 5 -9 10 12 20 24 -23 -1 -2 3 -4 -5 -6 -7 8 -10 -11}}
\end{figure}

\begin{figure}[h]
 \centerline{\includegraphics[width=\textwidth]{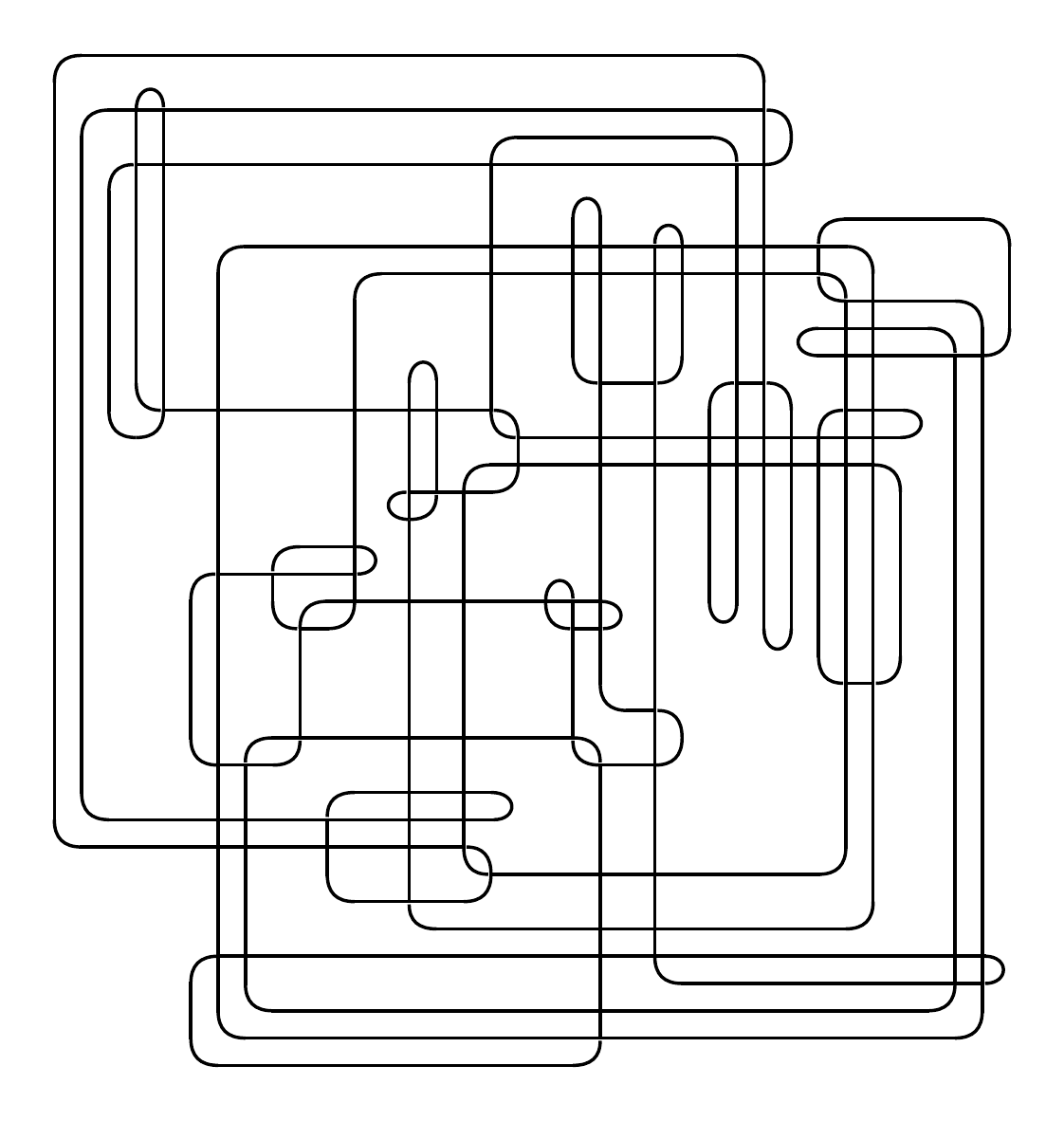}}
 \caption{$PZ_{120}$ {\em (\cite[Figure~12]{petronio2016algorithmic})}: \texttt{1 2 -120 -117 -116 -115 114 113 -112 -111 -110 32 -31 -30 29 40 -39 -38 37 78 77 -92 91 90 -89 88 87 68 -66 64 53 -51 48 45 -43 -41 38 34 -33 -37 41 42 -40 -36 35 39 -42 -44 46 49 -52 54 65 -67 69 116 108 99 -101 103 111 50 -49 -48 51 52 -50 -47 44 43 -45 -46 47 110 102 -100 98 12 11 80 -85 -91 92 86 -80 -79 83 89 -90 -84 79 9 10 -93 95 -98 100 101 -99 96 -94 97 -96 -95 93 94 -97 109 117 -74 71 -69 67 66 -68 70 -73 72 -70 -71 74 73 -72 81 82 -83 84 85 -86 75 76 33 -34 -35 36 25 -26 -27 28 -102 -103 -104 105 106 -107 -108 -109 6 -8 -10 -12 -14 17 20 -22 -24 27 31 -32 -28 24 23 -25 -29 30 26 -23 -21 19 16 -13 -11 -9 -7 5 -82 -88 57 -55 -77 -75 15 -16 -17 14 13 -15 -18 21 22 -20 -19 18 -76 -78 -56 58 -53 -54 118 -113 -105 104 112 -118 -119 115 107 -106 -114 119 -65 -64 -63 60 -58 56 55 -57 59 -62 61 -59 -60 63 62 -61 -87 -81 -1 3 -5 7 8 -6 4 -2 120 -4 -3}}
\end{figure}

\begin{figure}[h]
  \centerline{\includegraphics[width=\textwidth]{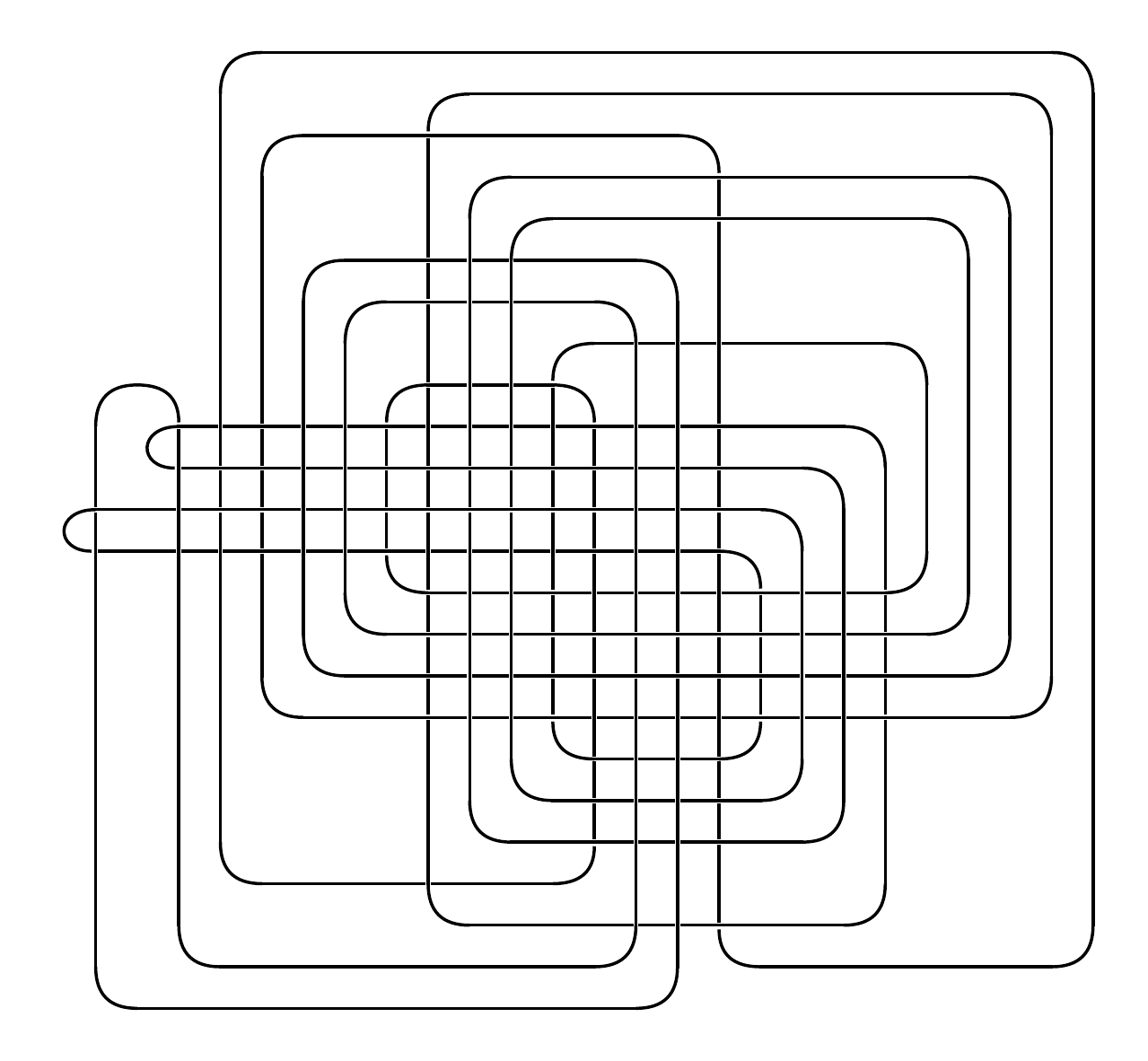}}
  \caption{$PZ_{138}$ {\em (\cite[Figure~27]{petronio2016algorithmic})}: \texttt{1 -11 -12 -13 -14 34 -46 -58 70 -86 -98 -110 -122 126 -130 -134 138 -111 -99 -87 -75 27 -28 -29 30 -31 -32 -33 -34 35 -36 -37 38 -138 -137 -136 -135 115 -103 -91 79 -63 -51 -39 -27 23 -19 -15 11 -38 -50 -62 -74 122 -121 -120 119 -118 -117 -116 -115 114 -113 -112 111 -1 2 -3 8 99 -100 -101 102 103 104 105 106 107 -108 -109 110 73 61 49 37 12 -16 -20 24 28 40 52 64 80 -92 -104 116 131 132 133 134 50 -49 -48 47 46 45 44 43 42 -41 -40 39 76 88 100 112 137 -133 -129 125 121 109 97 85 69 -57 -45 33 18 17 16 15 6 -4 -2 3 5 7 19 20 21 22 32 -44 -56 68 84 96 108 120 124 -128 -132 136 113 101 89 77 51 -52 -53 54 55 56 57 58 59 -60 -61 62 130 129 128 127 117 -105 -93 81 65 53 41 29 25 -21 -17 13 36 48 60 72 98 -97 -96 95 94 93 92 91 90 -89 -88 87 9 -5 4 -6 -7 -10 75 -76 -77 78 -79 -80 -81 -82 83 -84 -85 86 -71 -59 -47 -35 14 -18 -22 26 -30 -42 -54 -66 82 -94 -106 118 -123 -124 -125 -126 74 -73 -72 71 -70 -69 -68 -67 66 -65 -64 63 -78 -90 -102 -114 135 -131 -127 123 -119 -107 -95 -83 67 -55 -43 31 -26 -25 -24 -23 10 -9 -8}}
\end{figure}

\end{document}